\documentclass{article}
\usepackage[british]{babel}
\usepackage[utf8]{inputenc}
\usepackage{amsmath}
\usepackage{graphicx}
\usepackage[colorinlistoftodos]{todonotes}
\usepackage{hyperref}
\usepackage{multicol}
\usepackage{amsfonts}
\usepackage{lmodern}       
\usepackage{hyperref}
\usepackage[british]{babel}
\usepackage{enumitem}         
\usepackage{amsmath,amsbsy,amsfonts,amssymb}
\usepackage{amsthm}

\usepackage{times}
\usepackage{caption}
\usepackage{mathptmx} 
\usepackage{indentfirst}
\usepackage{float}
\usepackage{tikz}
\usepackage[square,numbers]{natbib}
\usepackage{hyperref}
\usepackage{subfigure}
\usepackage{colortbl}
\usepackage{xcolor}
\usepackage{mdframed}
\usepackage{asymptote}
\usetikzlibrary{shapes.geometric}
\newtheorem{theo}{Theorem}
\newtheorem{defi}{Definition}
\newtheorem{prop}{Proposition}
\newtheorem{coro}{Corollary}

\newtheorem{lema}{Lemma}
\usepackage{xcolor}
\usepackage{soul} % for \st
\usepackage{authblk}

\newtheoremstyle{meuestilo}  % nome do estilo
  {6pt plus 2pt minus 2pt}   % espaço antes
  {6pt plus 2pt minus 2pt}   % espaço depois
  {}  % fonte do corpo
  {1.5em}     % recuo (primeira linha)
  {\itshape}  % fonte do título
  {.}         % pontuação após o título
  { }         % espaço após o título
  {}          % especificação da "cabeça"

\theoremstyle{meuestilo}
\newtheorem{remark}{Remark}

\usepackage[many]{tcolorbox} % or \usepackage{tcolorbox} + \tcbuselibrary{skins}
\tcbuselibrary{skins}

\tcbset{
  myinlinebox/.style={
    enhanced,
    boxrule=0.5pt,
    colback=yellow!20,
    colframe=orange!80!black,
    sharp corners,
    boxsep=1pt,
    left=2pt,
    right=2pt,
    top=1pt,
    bottom=1pt,
    fontupper=\footnotesize,
  }
}

\newtcbox{\DD}{on line, myinlinebox}

\renewenvironment{proof}[1][Proof]%
{\noindent\textbf{#1.} }%
{\hfill\rule{0.7em}{0.7em}\par}
{\noindent\textbf{#1:} }%
{\hfill$\Box$\par}
\newmdenv[
    backgroundcolor=yellow!20, % Cor de fundo
    nobreak=true               % Não dividir o bloco entre páginas
]{paint}

\title{Topological Data Analysis in Finsler Spaces}
\author[1]{Rafael Cavalcanti}
\author[1]{Nelson Leal}
\author[2]{Danillo Barros de Souza}
\author[2]{Serafim Rodriguez}
\author[3]{Mathieu Desroches}

\affil[1]{Federal University of Pernambuco}
\affil[2]{Basque Center for Applied Mathematics}
\affil[3]{Inria, Montpellier, France}
\date{}

\begin{document}

\maketitle
\begin{abstract}
We introduce a new approach to Topological Data Analysis (TDA) based on Finsler metrics and we also generalize the classical concepts of Vietoris-Rips and Čech complexes within this framework. In particular, we propose a class of directionally dependent Finsler metrics and establish key results demonstrating their relevance to TDA. Moreover, we show that several Information Theoretic perspectives on TDA can often be recovered through the lens of Finsler geometry.

\end{abstract}

\tableofcontents

\section{Introduction}

Finsler geometry, named after the 1918 work of Paul Finsler \cite{finsler}, is a type of non-Euclidean geometry that presents a more general structure than Riemannian and Lorentzian geometries, for example, since the associated metric tensor depends on two variables, in addition to the position (point of the manifold), it also depends on the direction (tangent vector to the given position). Thus, it offers a richer, although more complex, structure that is already used in many contexts to model phenomena in different areas. Some those contexts in which Finsler geometry was used as a tool are: physics and biology in many scenarios \cite{ingarden, miron}, relativity \cite{relativity1, relativity2, relativity3}, mechanics \cite{mechanics1, mechanics2} and robotics \cite{robotics}. 

Topological data analysis (TDA) consists of using topological techniques to find structure in a data set \cite{tda1}, it started at the beginning of the century with pioneering works such as Edelsbrunner et al. in 2002 \cite{tda2} e Zomorodian e Carlsson in 2005 \cite{tda3}. Persistent homology of a filtration of a family of simplicial complex is the most commonly tool used topological method to identify structure in a data set, although there are studies that use other methods from a broader perspective \cite{tda4}. In recent years, TDA has been used for data analysis in many contexts, some of which are: image processing \cite{images1, images2, images3}, time series analysis \cite{timeseries1, timeseries2, timeseries3} and neuroscience \cite{neuro1, neuro2, neuro3}.

Data can be understood as a set of points with numerical entries, each originating from a selected feature. These points can therefore be visualized in a subset $\Omega$ of the Euclidean space $\mathbb{R}^n$. With this perspective, we can consider a notion of distance in $\Omega$ and infer a notion of proximity between the data. The topology of the data is obtained through homological persistence diagram of a filtration of discrete structures called simplicial complexes. Under certain conditions, these structures capture the same type of topological information as the dataset, and the most commonly used simplicial complexes for this purpose are the Čech complex and the Vietoris–Rips complex.

The Čech complex is constructed by considering higher-order interactions among the points, but this structure can be approximated by the Rips complex, which is built from pairwise notions of distance between points. Constructing these structures by considering only distances has proven effective in many TDA applications, but, as with any attempt to establish a signature for a dataset, some information may be lost in the process \cite{tda1}. This loss is particularly likely when the metric fails to capture the underlying relational geometry or the contextual significance among data points. In particular, relying solely on pairwise distances without incorporating additional information about the points may obscure meaningful higher-order interactions or anisotropic behaviors present in the data. This motivates the search for stronger, more expressive frameworks that extend or refine the role of distance in TDA and that can capture interactions more consistently with the data. The figure below is a good sketch of how weird, but still convex, might the balls built by a Finsler metric be.
\begin{figure}[H]
    \centering
    \includegraphics[width=1.0\textwidth]{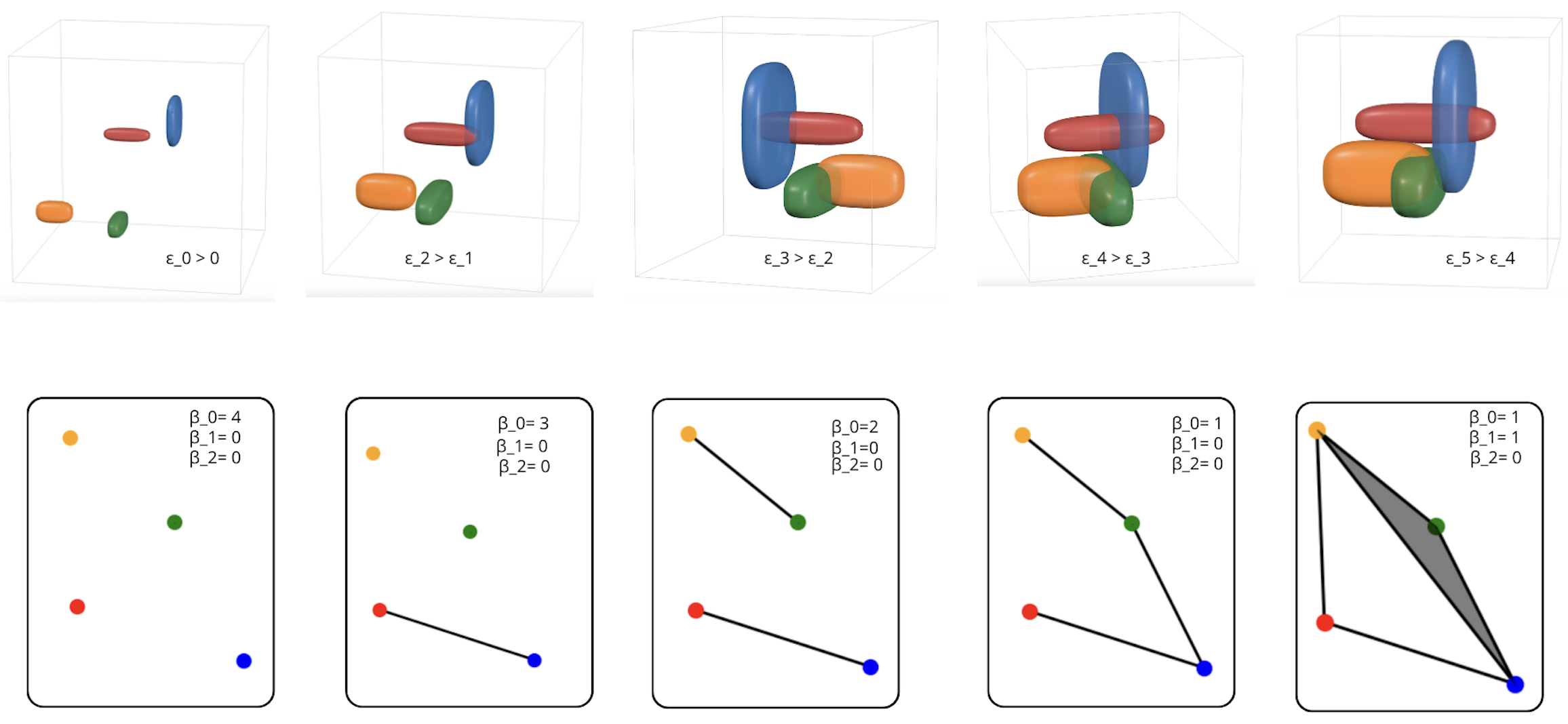}
    \caption{On the top part of the figure: 4 random points in $\mathbb{R}^3$, with balls of same radius centered at these points, built by a Finsler metric. On the bottom part: the Cech complexes for each $\epsilon$ in a filtration.}
    \label{fig:my_image}
\end{figure}
\vspace{-0.5em} 
Several works in population dynamics have been published from Antonelli's perspective \cite{verdinho, vh1, vh2}, taking in consideration the benefit of a Finsler metric be a function of two variables. In Atonelli's approach, when $F(p,v)$ is a Finsler metric, $p$ plays the role of being the production of species and $v$ being the population density. Looking at this approach, arises the idea of using a Finsler structure to study TDA based on two important characteristics of the data: the distances between the data points and the intrinsic information about the points at which we calculate the distance, basically, we extract the point and its direction.

In this work, we address the topological analysis of data by investigating distance functions in the context of a Finsler metric or quasi-metric space and analyzing their impact on the construction of simplicial complexes that are essential for obtaining the homological persistence diagram. We stopped investigating the already consolidated Euclidean distance and moved on to the analysis of more general distance functions. We began by inserting into a Finsler context distances known in the literature that, even though they are not Euclidean, are still isometric to it. Our main results arise in a context that allows us to see the data points also as tangent vectors and the Finsler metric is used to generate distance functions that consider the contextual significance between the points, thus, from the information present in the data points, we can begin to capture significant interactions. With the proposed non-standard distances, we consider generalized Cech and Vietoris-Rips complexes, and propose a new type of simplicial complex based on the pairwise relationship between data points that naturally carries more information about the points than the distance between them. This new concept of simplicial complex cannot search for such significant high-order interactions as the Cech complex, but it certainly can search for more interactions between pairs than the Vietoris-Rips complex. Finally, considering the computational complexity of constructing the generalized Cech complex, we present a generalization of an already known inclusion involving the Vietoris-Rips complex that allows us to obtain a good approximation for the Cech complex from lower and upper bounds.

The structure of this paper is organized as follows. In Section \ref{section_finsler}, we present general tools from Finsler geometry that will be used throughout our results and approach. In Section \ref{sec_iso_convex}, we engage with results from a recent work on TDA and information spaces, reinterpreted through the lens of Finsler geometry; in addition, we propose an improvement to one of their results.  Section \ref{section_tda} is devoted to a brief and broad overview of Topological Data Analysis (TDA), providing the necessary background. Section \ref{section_fisnler_tda} introduces new ideas for developing TDA in Finsler spaces using a distinct approach, where we formulate new definitions and prove novel results. Finally, in Section \ref{section_stability}, we adapt concepts from the classical theory of persistence stability to the Finslerian setting, extending known results in TDA to this more general geometric framework.

%-----------------------------------------------------
\section{Finsler Structure}\label{section_finsler}
%-----------------------------------------------------

In this section we will assume some notion of general structures such as manifolds, whose we will simply face it as a topological space that is locally Euclidean, in the sense that for each point on this space there is a diffeomorphism from a neighborhood of this point to an open set of $\mathbb{R}^n$. Beyond that, if $M$ is a manifold, we denote by $TM$ the disjoint union of all tangent space of $M$ and the null section $\mathbb{O}$ of $TM$ will be the space built as the union of each null vector of all of tangent space of $M$. 

\begin{defi}\label{defiespfinsler}
    Let $M$ be a manifold and $\mathbb{O}$ the null section of $TM$. A continuous function $F : TM \rightarrow \mathbb{R}$ is a \textbf{Finsler metric} if it satisfies
\begin{enumerate}[label=(\roman*)]
    \item $F$ is of $C^\infty$ class on $TM \backslash \mathbb{O}$;
    \item $F$ is positive definite on $TM \backslash \mathbb{O}$;
    \item $F$ is p-homogeneous of degree $1$ in $v$, that is, $F(p, rv) = r F(p,v),~\forall r>0$;
    \item The matrix with coefficients $g_{ij}(p,v) = \dfrac{1}{2} \dfrac{\partial^2 F^2}{\partial v^i \partial v^j}$ is positive definite on $TM \backslash \mathbb{O}$. 
\end{enumerate}
    If $F$ is a Finsler metric, the pair $(M,F)$ is called a \textbf{Finsler space}. 
\end{defi}

%A Finsler metric is defined on tangent bundle and then it depends of two variables that in general are view as position and velocities or other measures as: populational density and production measure \cite{}\textcolor{red}{colocar citações metricas de finsler}. 

% \DD{Explain the details of considering the point of application $p$ and the tangent vector $v$ and their importance on Finsler metrics.}

There is also a Finsler version of the Euclidean metric described by the function $E:TM \rightarrow \mathbb{R}$, given by 
\[
E(p,v) = \left( \sum_{i=1}^n (v^i)^2\right)^{\frac{1}{2}},
\]
this is called \textbf{Euclidean Finsler metric}. We have that $E(p,rv)=|r|E(p,v), \forall r \in \mathbb{R}$, the metrics that have this property are called absolutely homogeneous in $v$. Beyond that $E^2$ has the property of decomposing in terms of $v^i$ only, where $v^i$ if the $i$-th coordinate of vector $v$. In general, if $f:U\times V\rightarrow\mathbb{R}$ takes the form $f(u,v)=\sum_{i=1}^nf_i(u,v_i)$, we say that $f$ is \textbf{decomposable} in the variable $v$, or we simply say $f$ is decomposable. 

\begin{defi}\label{defi_comprimento}
    The \textbf{absolute length} $\| v \|_{p}$ of a tangent vector $v=v^i\frac{\partial}{\partial p^i}$ in $T_pM$ is defined as the value of $F$, that is, $\| v \|_p = F(p,v)$. 
\end{defi}

\begin{remark} \label{remark1}
By Definition \ref{defiespfinsler}, $F$ is positive definite, the triangular inequality for $F$ is valid because $g_{ij}$ is positive definite then $F$ is a convex function. Note that we do not require $F$ to be absolutely homogeneous, i.e., we may be $F(p,-v) \neq F(p,v)$, that is why $F$ is a quasi-norm in $T_pM.$
\end{remark}

\begin{defi}
    The \textbf{length of a curve} $\gamma: [a,b] \rightarrow M$ with endpoints $\gamma(a) = x$ and $\gamma(b) = y$ is defined as the value 
\begin{equation}\label{comprimentodecurva}
    L(\gamma) = \int_a^b F \left( \gamma,\dot{\gamma}\right ) dt.
\end{equation}
    The \textbf{distance between two points} $x$ and $y$ in a Finsler space is defined as
\begin{equation}\label{distancia}
    d_F(x,y) = \inf_{\gamma} ~ L(\gamma),
\end{equation}
    where $\gamma$ is any curve connecting the points $x$ and $y$.
\end{defi}

The distance (\ref{distancia}) is a natural distance of a Finsler space, as in any curve space. From Remark \ref{remark1} we can have $d_F(x,y) \neq d_F(y,x)$, that is, $d_F$ is a quasi-metric. An alternative to this problem is to take $$d(x,y) = \frac{d_F(x,y)+d_F(y,x)}{2}$$ and the tuple $(M,F,d)$ will be a metric space. But, if $F$ is absolutely homogeneous of degree 1 in $v$ (as the Euclidian Finsler metric, for example), then we do not have this problem and $d_F$ is a metric \cite{tamassy}, so that, $(M,F,d_F)$ is a metric space.  

\begin{remark}
From here on, we will denote points in $M$ in two ways: when calculating the value of $F$, we will use $p$ and $q$; when calculating the distance between two points, we will use
$x$ and $y$. This distinction will become relevant in the following sections.
\end{remark}

\begin{defi}
A family of curves $\{\gamma_s\}$ for $s \in [-\epsilon , \epsilon]$, all defined in a, interval $[a,b]$ with $\gamma_s(a)=p$ and $\gamma_s(b)=q$, is smooth if the map
$$
\begin{array}{cccc}
& (-\epsilon, \epsilon) \times [a, b] 
& \longrightarrow
& M \\
& (s,t)
& \longmapsto
& \gamma_s(t)
\end{array}
$$
is a smooth map. A curve $\gamma$ is an extremal of $L(\gamma)$ if for every smooth family $\left\{\gamma_s\right\}$ with $\gamma = \gamma_0$ we have $$\displaystyle\left.\frac{d L\left[\gamma_\tau\right]}{d \tau}\right|_{\tau=0}=0$$.
\end{defi}

If a curve $\gamma$ is an extremal of $L(\gamma)$, then the coordinates $\gamma^i(t) = (p^i(t), \cdots, p^n(t))$ satisfies the Euler-Lagrange equations
\begin{equation}\label{eulerlagrange}
    \frac{d}{dt}\left(\frac{\partial F^2}{\partial v^i}\right) = \frac{\partial F^2}{\partial p^i}, \quad \mbox{ for} \quad v^i = \frac{dp^i}{dt}.
\end{equation} 
Note that the solutions of \eqref{eulerlagrange} describing the curve that perform the distance on $(M,F)$. In truth, these are the equations of the geodesics of $M$, for the metric $F$. 

For the Euclidean Finsler metric, the curve $\gamma$ that satisfies the solutions of Euler-Lagrange equations are $\gamma^i(t) = a_i~t + b_i$ for constants $a_i$ and $b_i$. For two points $x,y \in M$ with $\gamma(0) = x^i$ and $\gamma(1) = y^i$, we have $a_i = y^i - x^i$ and $b^i = x^i$, then
\begin{equation}\label{euclidiandist}
d_E(x,y) = \int_0^1 E(a_i~t + b_i, a_i)~dt = \int_0^1 E(t(y^i - x^i) + x^i, y^i - x^i)~dt  = \left(   \sum_{i=1}^n (y^i - x^i)^2\right)^{\frac{1}{2}},
\end{equation}
and the infimum of (\ref{distancia}) is omitted because we are considering geodesic curves. 
    
We will investigate some metrics isometric to the Euclidean metric, so that we need to give a suitable definition for when two Finsler metrics are isometric. For us, isometry in Finsler spaces are the same used in \cite{venegas}, and with this well-defined concept we will show later that the isometry on the distances induced by Finsler metric can be rescued by the isometry on the Finsler structures.
    
\begin{defi}
    A mapping $\varphi:(M,F) \rightarrow (N,G)$ between two Finsler spaces is said to be a \textbf{Finsler isometry} if it is a diffeomorphism which preserves the Finsler structure, that is, for every $p \in M$ and every $v \in T_pM$, we have 
\begin{equation}
    F(p,v) = G(\varphi(p), d\varphi_p(v)). 
\end{equation}
\end{defi}

In \cite{isometrygroup} the authors give a detailed discussion about isometries in Finsler spaces, relating them with the induced distance on the base manifold. The following result is the bridge between the isometries in different contexts.
    
\begin{theo}\label{isometrydistance}
    Let $(M,F)$ and $(N,G)$ be connected Finsler spaces and consider a bijection $\varphi:(M,F) \rightarrow (N,G)$. Then $\varphi$ is a Finsler isometry if and only if it is an isometry for the corresponding Finsler distance, that is, for any $p,q \in M$, we have 
\begin{equation}
    d_F(x,y) = d_G(\varphi(x), \varphi(y)). 
\end{equation}
\end{theo} 

Let $E$ be the Euclidean Finsler metric. If $\varphi: M \rightarrow M$ is a Finsler isometry between $(M, F)$ and $(M, E)$, then using the Theorem \ref{isometrydistance}, for all $x,y \in M$, we have
\begin{equation}\label{distanceF_isometry}
    d_F(x,y) = d_E(\varphi(x), \varphi(y)) = \left( \sum_{i=1}^n \big(\varphi(y^i) - \varphi(x^i)\big)^2\right)^{\frac{1}{2}}.
\end{equation}

Isometries between $F$ and $E$ simplify very much the computation of the induced distance on $(M,F)$, since $d_F$ will be given by \eqref{distanceF_isometry}. 

In \cite{fisher}, it was shown that, in the context of information geometry, certain distances derived from divergences are isometric to the Euclidean metric. In this work, we present results that establish analogous isometries for such distances, now obtained by means of Finsler metrics.

\begin{defi}
    A Finsler space $(M,F)$ is called \textbf{locally Minkowski space} if $M$ is covered by coordinates neighborhoods $(\bar{p}^i)$ such that $F$ does not depend on $(\bar{p}^i)$. Such coordinates $(\bar{p}^i)$ are called \textbf{adapted}.
\end{defi}

\begin{lema}\label{lemma1}
    Let $\Omega\subset \mathbb{R}^n$ be an open connected set and $E: T\Omega \rightarrow \mathbb{R}$ be the Euclidean Finsler metric. A Finsler metric $F: T\Omega \rightarrow \mathbb{R}$ is isometric to $E$ if and only if $F$ is locally Minkowski and $F^2$ is decomposable.
\end{lema}

\begin{proof}
    Isometry implies existence of a diffeomorphism $\varphi:\Omega \rightarrow \Omega$ such that 
    \[
    F^2(p,v) = E^2(\varphi(p), d\varphi_p(v)),
    \]
    for every $p\in M$ and every $v\in T_pM$. Let $\phi = \varphi^{-1}$ and consider the coordinate system $\bar{p}^i = (\phi(p^i))$. Thus, 
    \[ 
    F^2(\bar{p}, \bar{v}) = F^2(\phi(p),d_{p} \phi(v)) = E^2\left(\varphi(\phi(p)), d_{\phi(p^i)} \varphi \cdot d_p \phi (v)\right) = \displaystyle\sum_{i=1}^{n}(v^i)^2, 
    \]
    that does not depend on $\bar{p}^i$. 

    Reciprocally, if $F$ is Minkowski, then there exists  $(\bar{p}^i)$ such that $F(\bar{p}, \bar{v})=F(\bar{v})$. On the other hand, $F^2$ decomposable may be written such as $F^2(\bar{v})=\sum_{i=1}^{n}F_i^2(\bar{v}^i)$. Setting $p^i = \varphi_i (\bar{p}^i)$, then $v^i=\varphi'_i(\bar{p}^i)\bar{v}^i$. Thus,
    \[
    E^2(p,v)=\sum_{i=1}^n (v^i)^2 = \sum_{i=1}^n (\varphi_i'(\bar{p}^i))^2(\bar{v}^i)^2.
    \]
    Using the equality above, the following ODE is required to obtain an isometry
    \[
    \left(\varphi_i^{\prime}\left(\bar{p}^i\right)\right)^2\left(\bar{v}^i\right)^2 = F_i^2\left(\bar{v}^i\right).
    \]
    Since $F^2$ is strictly convex, the function $\displaystyle\psi(\bar{p}^i)=\int_0^{\bar{p}^i}\frac{F_i(\xi^i)}{\xi^i}d\xi^i + c_0$ is solution of the previous ODE.
\end{proof}

It is not so easy to say when $d_F$ and $d_E$ are isometric, so then, provided $F$ and $F^2$ satisfy the hypothesis of Lemma \ref{lemma1}, it is sufficient to get the isometry.

%-----------------------------------------------------
\section{Isometry, convexity and contractibility}\label{sec_iso_convex}
%-----------------------------------------------------

Let \( X \subset \mathbb{R}^n \) be an open connected set. We say that \( X \) is \textbf{convex} if, for all \( x, y \in X \) and for every \( \lambda \in [0,1] \), the point on the line segment given by \( (1-\lambda)x + \lambda y \) lies entirely within \( X \).

A basic but important property of convex sets is that their intersection remains convex. More precisely:

\begin{prop}\label{int_convex}
    If \( X \) and \( Y \) are convex sets, then their intersection \( X \cap Y \) is also convex.
\end{prop}

Convexity also has implications in topology. A topological space \( X \) is said to be \textbf{contractible} if it is homotopy equivalent to a single point, that is, there exists a point \( p \in X \) such that \( X \simeq \{p\} \).

We now state a useful result connecting convexity and contractibility:

\begin{prop}\label{convex_to_contra}
    If \( X \) is a convex set, then \( X \) is contractible.
\end{prop}

\begin{remark}
This is the important result to guarantee the modified version of the Nerve Theorem, commented in the Remark \ref{remark2}. 
\end{remark}
 
Let $\mathbb{R}_{>0}^n$ be a subset on the strictly positive orthant of $\mathbb{R}^n$, that is, all coordinates are strictly positive $\mathbb{R}_{>0}^n = \{(p^1,\ldots, p^n) \in \mathbb{R}^n;p^i>0, \forall i  \}$. The function  $F:\mathbb{R}_{>0}^n \times \mathbb{R}^n \rightarrow \mathbb{R}$ defined by
\begin{equation}\label{lnmetric1}
    F(p,v) = \left( \sum_{i=1}^{n}\left( \frac{v^i}{p^i}\right)^2 \right)^\frac{1}{2}, 
\end{equation}
is a Finsler metric. Note that $F$ is absolute homogeneous of degree 1 in $y^i$, then $\mathbb{R}^n$ is a metric space with the distance $d_F$ induced by $F$ as in \eqref{distancia}. Let $\varphi: \mathbb{R}_{>0}^n \rightarrow \mathbb{R}_{>0}^n$ with $\varphi_i(p^i) = \ln{(p^i)}$ and note that $\varphi_i'(p^i) v^i = \frac{v^i}{p^i}$, then
\[
E(\varphi(p), d\varphi_p(v)) = \left( \sum_{i=1}^n \big(\varphi_i'(p^i) v^i\big)^2\right)^{\frac{1}{2}}  = \left( \sum_{i=1}^{n}\left( \frac{v^i}{p^i}\right)^2\right)^\frac{1}{2} = F(p,v),
\]
i.e., $\varphi$ is a Finsler isometry between $F$ and $E$, thus by \eqref{distanceF_isometry},
\begin{equation}\label{distanceF_ln}
    d_F(x,y) = \left( \sum_{i=1}^n \big(\ln(y^i) - \ln(x^i)\big)^2\right)^{\frac{1}{2}}, 
\end{equation}
this is the distance for \textit{Burg information metric} as described in the paper \cite{fisher}. 

\begin{remark}
The identity (\ref{distanceF_ln}) can be obtained using the Lemma \ref{lemma1}, just observe that $F$ given by \eqref{lnmetric1} is such that $F^2$ is decomposable and take the adapted coordinates given by $\bar{p}^i = \ln({p^i})$ for show that $F$ is locally Minkowski. 
\end{remark}

In \cite{fisher}, they investigate the balls $B_r^D(x)$, centered at $\textbf{1}=(1,\ldots,1)$ of the \textit{Itakura-Saito divergence}, $D(x \| y)=\sum_{i=1}^n\left[\ln \frac{y_i}{x_i}+\frac{x_i}{y_i}-1\right]$, and then conclude that those balls are convex for radius $r^2 \leq \ln 2-\frac{1}{2}$ and nonconvex if $r^2>2 \ln 2-1$. Although, we realized that built balls using the distance \eqref{distanceF_ln} have large radius of convexity in comparison with the divergence in \cite{itakura}.
    
\begin{theo}
    Let $\Omega \subseteq \mathbb{R}^n$ be an open connected domain, such that $\partial \Omega$ is an orientable surface of codimension 1. Then $\Omega$ is convex if and only if the second fundamental form of $\partial \Omega$ is positive semi-definite.
\end{theo}

\begin{coro}\label{sinal_curvatura}
    Let $\Omega \subseteq \mathbb{R}^n$ be an open connected domain. Then $\Omega$ is convex if and only if the Gaussian curvature of $ \partial \Omega$ does not change its signal.
\end{coro}

\begin{proof}
    Since positive semi-definite operators have non-negative eigenvalues, and the Gaussian curvature of $\partial \Omega$ is product of these eigenvalues, the result follows.
\end{proof}

\begin{prop}
    Let $d_F$ be as in \eqref{distanceF_ln}. The ball $B_{d_F}(\mathbf{1},\rho)$ is convex for $\rho\leq1$.
\end{prop}

\begin{proof}
    We could follow the canonical steps to prove that a general region is convex, using the line segment definition, but we opted for a more elegant proof using geometric concepts. Consider the function   $f\left(x^1, x^2, \ldots, x^n\right)=\sum_{i=1}^n \ln ^2\left(x^i\right)-\rho^2$. The Gaussian curvature of the level-set surface given by $f\left(x^1, x^2, \ldots, x^n\right)=0$ is      $K=\frac{\operatorname{det}\left(H_f\right)}{\|\nabla f\|^4}$, where $H_f$ is the hessian of $f$. Omitting the computations, we obtain 
    \[
    K=\frac{\prod_{i=1}^n 2\left(1-\ln \left(x_i\right)\right)}{\left(\sum_{i=1}^n 4 \ln ^2\left(x_i\right)\right)^2} \cdot \prod_{i=1}^n \frac{1}{x_i^2},
    \]
    which has non-negative values for $x^i\leq e$. Let $x^i=1, i\neq j$, we get $f\left(1, \ldots,x^j, \ldots, 1\right)= \ln^2(x^j)-\rho^2$.  If $\rho$ vary on $(0,1]$, then $\ln ^2\left(x^j\right)-\rho^2\leq 0$ for $x^j \in [e^{-1},e]$. This easily implies that $B_d(\mathbf{1}, \rho)$ is convex if $\rho\leq 1$.
\end{proof}

\begin{prop}
    Let $d_F$ be as in \eqref{distanceF_ln}. The balls generated by $d_F(x,y)$ have the property that every finite intersection is either empty or contractible.
\end{prop}
\begin{proof}
    Without loss of generality lets consider $\mathbf{1}=(1, \ldots, 1)$ and $a = (a^1, \ldots, a^n)$. For the purpose of this proof, consider the maps 
    \[
    \log : \mathbb{R}_{>0}^n \rightarrow \mathbb{R}^n, \quad x=\left(x_1, x_2, \ldots, x_n\right) \mapsto\left(\ln (x_1), \ln (x_2), \ldots, \ln (x_n)\right)
    \]
    with its respective inverse map given by 
    \[
    \exp : \mathbb{R}^n \rightarrow \mathbb{R}_{>0}^n, \quad y=\left(y_1, y_2, \ldots, y_n\right) \mapsto\left(\exp \left(y_1\right), \exp \left(y_2\right), \ldots, \exp \left(y_n\right)\right).
    \]
    
    The image of $B_{d_F}(\mathbf{1}, \rho)$ under $\log$ map is a Euclidean ball described as the set
    \[
    \log \left(B_{d_F}(\mathbf{1}, \rho)\right)=\left\{y=\left(\ln (x_1), \ln (x_2), \ldots, \ln (x_n)\right): \sum_{i=1}^n y_i^2 \leq \rho\right\},
    \]
    and hence 
    \[
    \log \left(B_{d_F}(\mathbf{1}, \rho) \cap B_{d_F}(a, \rho)\right)=\log \left(B_{d_F}(\mathbf{1}, \rho)\right) \cap \log \left(B_{d_F}(a, \rho)\right)=\mathbb{B}(0, \rho) \cap \mathbb{B}(\log a, \rho),
    \]
    where $\mathbb{B}(p, \rho)$ is the Euclidean ball of radius $\rho$ centered at $p$. Since non empty intersection of Euclidean balls is a convex set, and convex sets are contractible,  the intersection $B_{d_F}(\mathbf{1},\rho)\cap B_{d_F}(a, \rho)$ is contractible once the exponential map preserves contractibility. For every non empty finite intersection of balls generated by $d_F$, one get contractibility in the Euclidean space and then by the inverse homeomorphism of $\log$ one obtain a contractible space as the intersection of balls given by $d_F$.
\end{proof}

\begin{figure}[H]
    \centering

    \begin{minipage}[t]{0.19\textwidth}
        \centering
        \includegraphics[width=\linewidth]{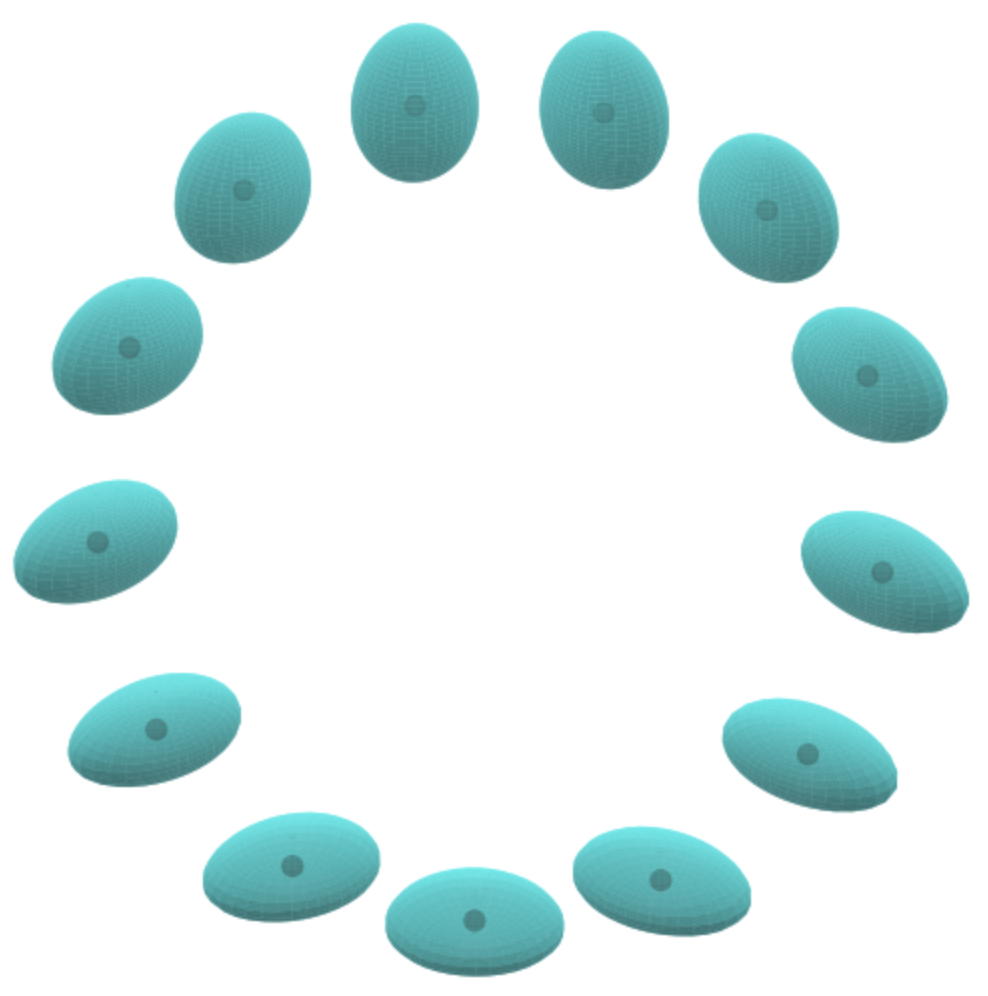}
        \caption*{ $\rho$ = 0.08}
    \end{minipage}
    \hfill
    \begin{minipage}[t]{0.19\textwidth}
        \centering
        \includegraphics[width=\linewidth]{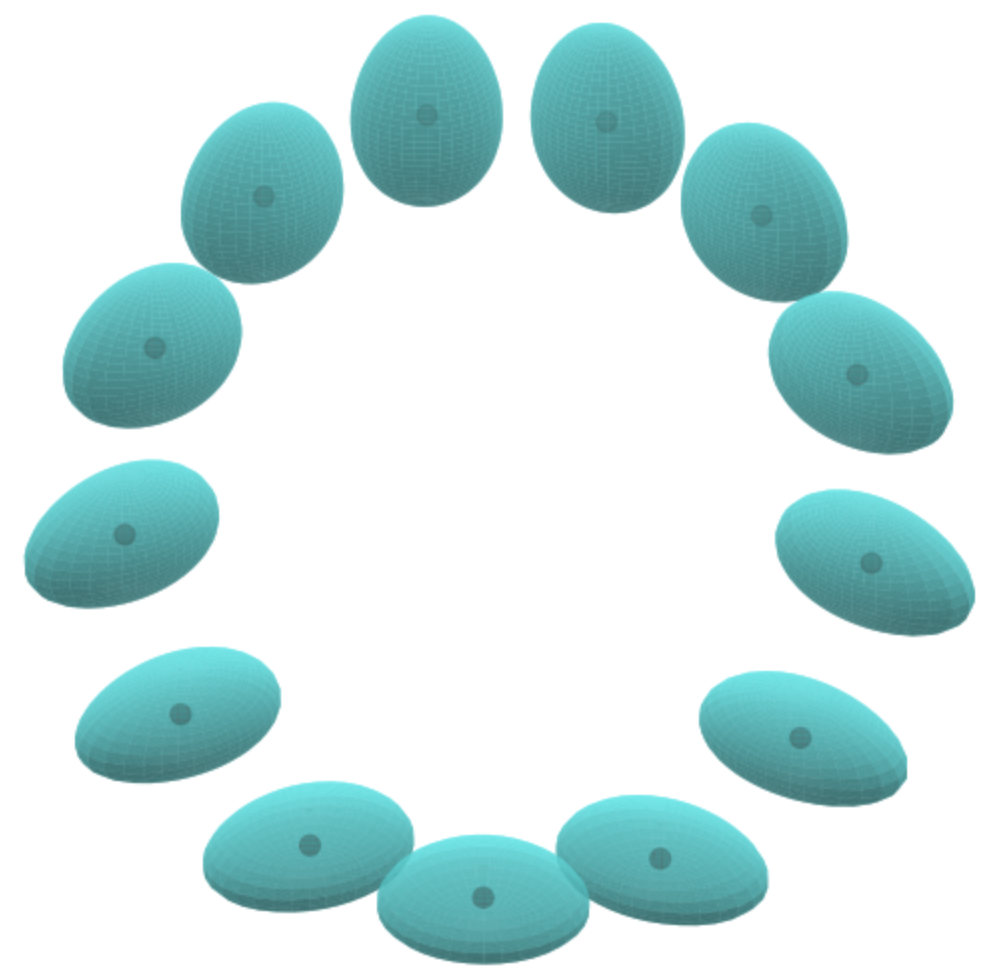}
        \caption*{ $\rho$ = 0.1}
    \end{minipage}
    \hfill
    \begin{minipage}[t]{0.19\textwidth}
        \centering
        \includegraphics[width=\linewidth]{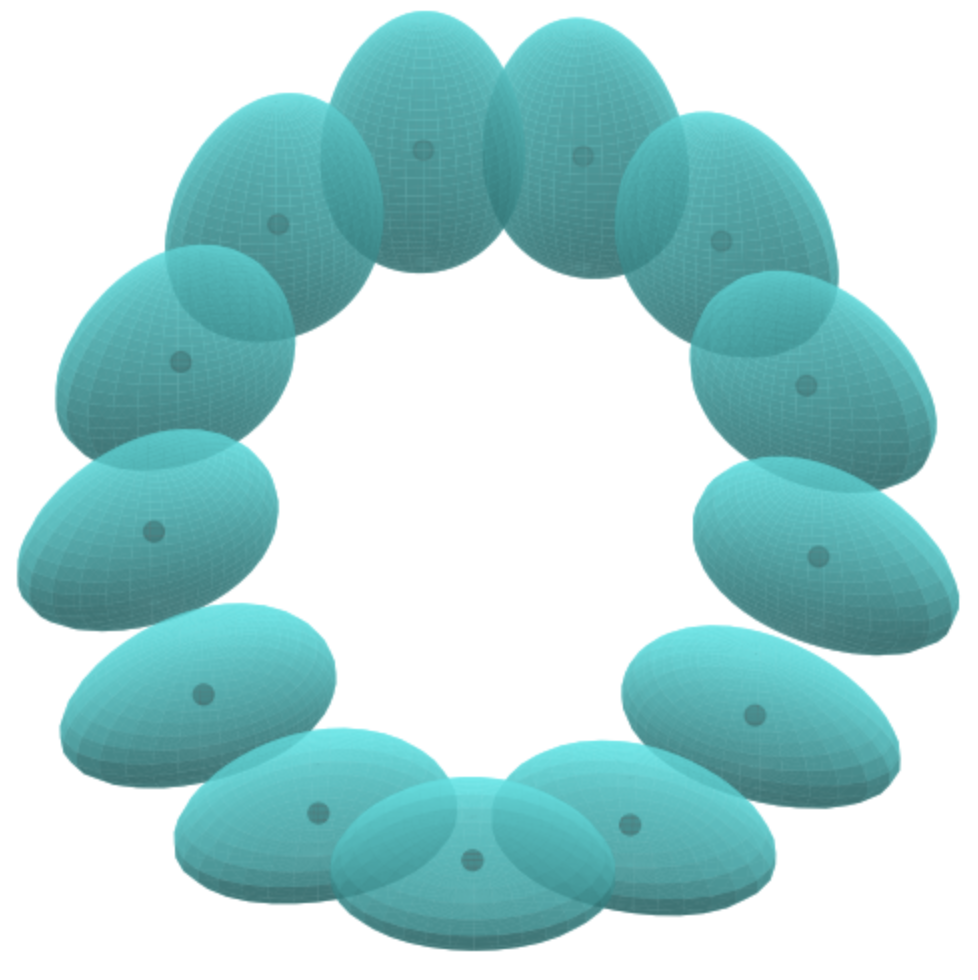}
        \caption*{ $\rho$ = 0.15}
    \end{minipage}
    \hfill
    \begin{minipage}[t]{0.19\textwidth}
        \centering
        \includegraphics[width=\linewidth]{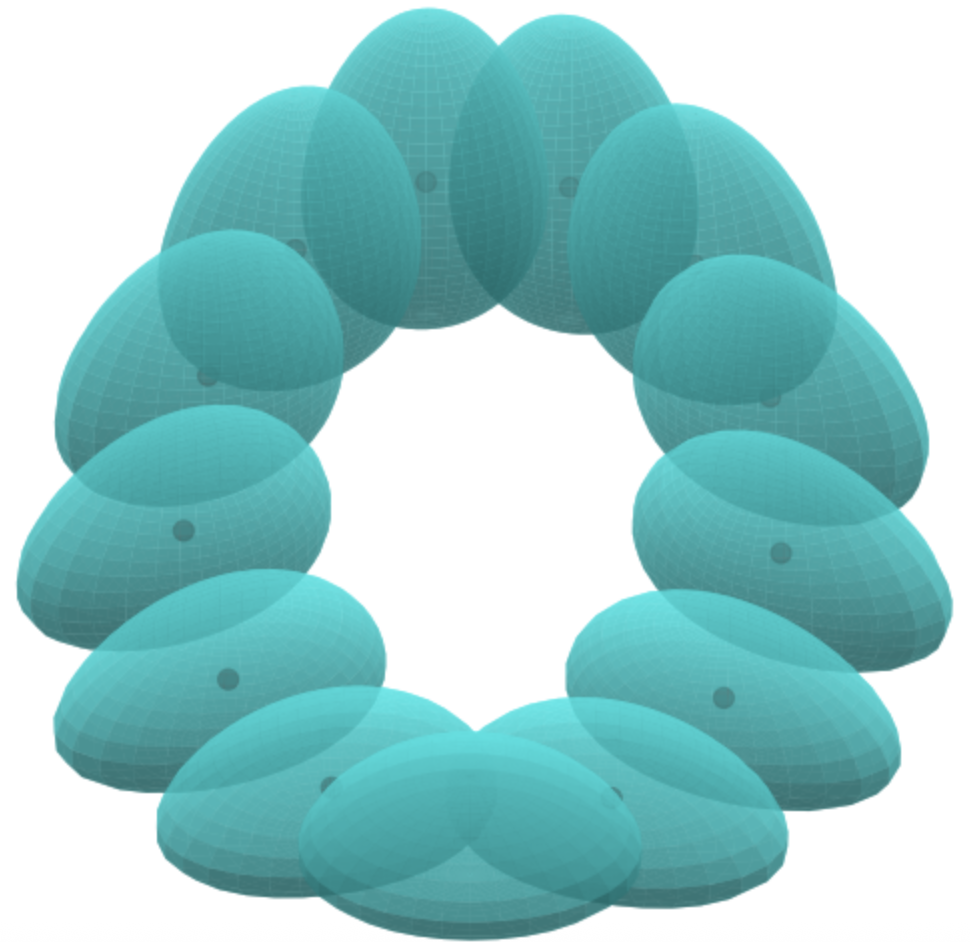}
        \caption*{ $\rho$ = 0.2}
    \end{minipage}

    \caption{The shape of the balls given by distance \eqref{distanceF_ln}.}
    \label{fig:log-comparison1}
\end{figure}

In Figure \ref{fig:log-comparison1} we use the distance \eqref{distanceF_ln} and a point set lying on a circle for illustrates the evolution of balls centered at each point for different radius. 

In \cite{fisher} the authors used \textit{Antonelli}'s isometry \cite{antonelli_isometry} to show that the distance derived from the \textit{Shannon entropy} is isometric to the Euclidean distance. Once again we are in position of formulating this discussion in terms of a Finsler metric. 

Let $F$ be another Finsler metric, but now defined by
\begin{equation}\label{metrica_raiz}
    F(p,v) = \left( \displaystyle\sum_{i=1}^{n}\left( \frac{v^i}{2\sqrt{p^i}}\right)^2 \right)^\frac{1}{2}. 
\end{equation}

The metric $F$ given by \eqref{metrica_raiz} is also a Finsler metric absolute homogeneous of degree 1 in $y^i$ and $\mathbb{R}_{>0}^n$ endowed with the distance $d_F$ induced by \eqref{metrica_raiz} is a metric space. Let $\varphi: \mathbb{R}_{>0}^n \rightarrow \mathbb{R}_{>0}^n$ with $\varphi_i(p^i) = \sqrt{p^i}$ and note that $\varphi$ is a Finsler isometry between $F$ and $E$, then again by \eqref{distanceF_isometry},
\begin{equation}\label{distanceF}
    d_F(x,y) = \left( \sum_{i=1}^n \Big(\sqrt{y^i} - \sqrt{x^i}\Big)^2\right)^{\frac{1}{2}}, 
\end{equation}
this is the distance for \textit{Fisher information metric} \cite{fisher}, and the balls associated to this distance are convex for any radius $\rho$ as they proved in the paper. As convexity implies contractibility, one can also says that the balls generated by the distance above  has non empty finite intersection contractible. 

\begin{figure}[H]
    \centering

    \begin{minipage}[t]{0.19\textwidth}
        \centering
        \includegraphics[width=\linewidth]{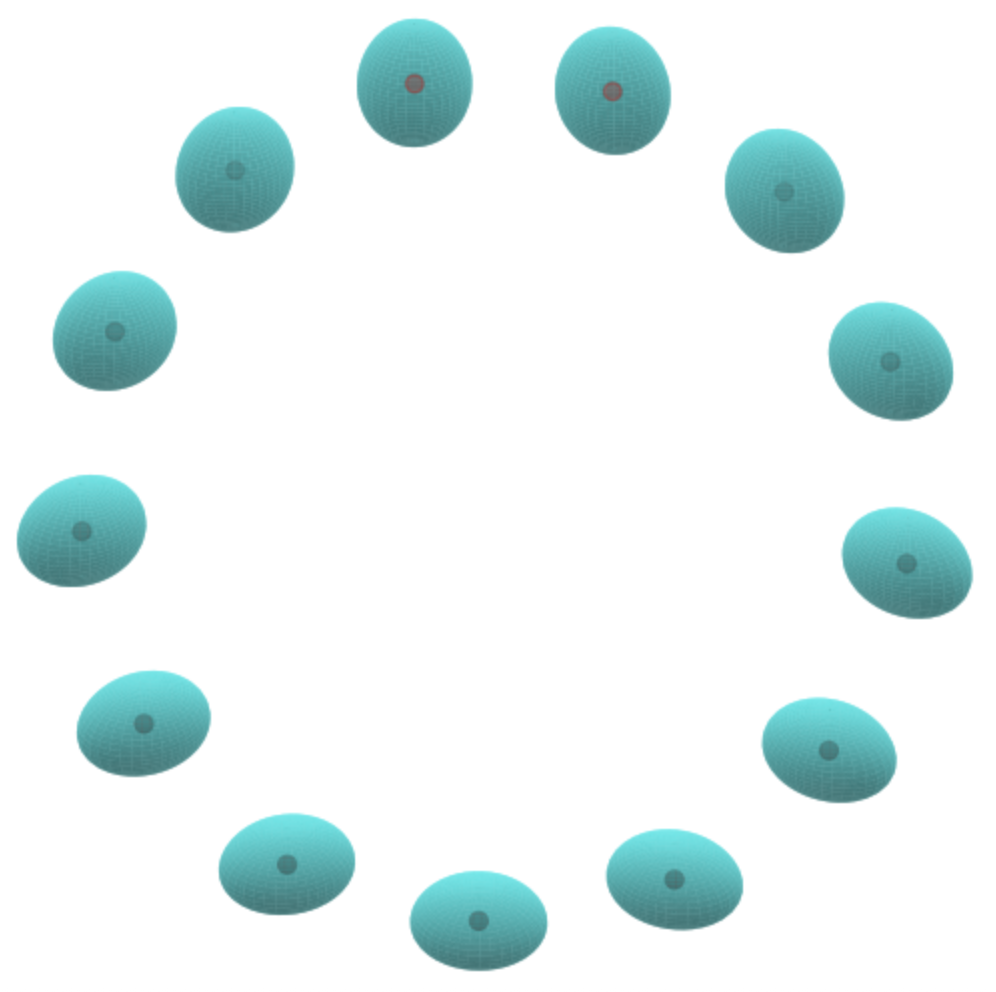}
        \caption*{$\rho$ = 0.12}
    \end{minipage}
    \hfill
    \begin{minipage}[t]{0.19\textwidth}
        \centering
        \includegraphics[width=\linewidth]{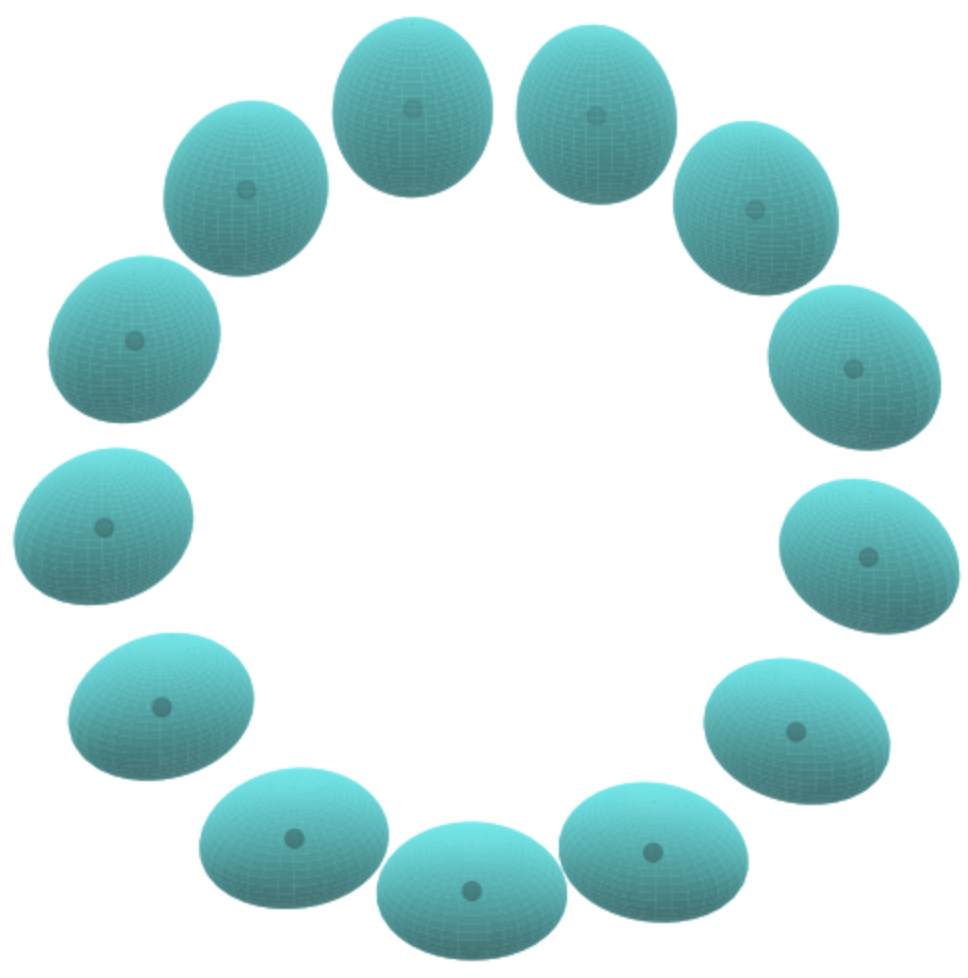}
        \caption*{ $\rho$ = 0.18}
    \end{minipage}
    \hfill
    \begin{minipage}[t]{0.19\textwidth}
        \centering
        \includegraphics[width=\linewidth]{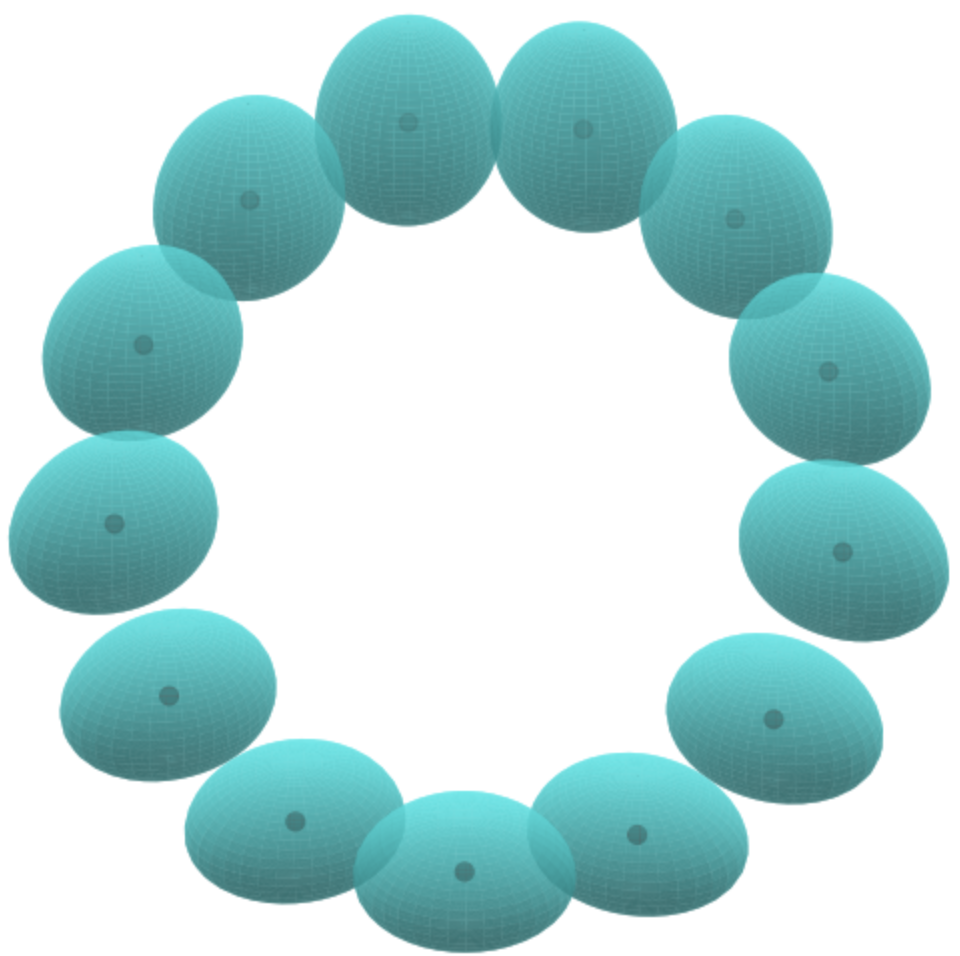}
        \caption*{ $\rho$ = 0.22}
    \end{minipage}
    \hfill
    \begin{minipage}[t]{0.19\textwidth}
        \centering
        \includegraphics[width=\linewidth]{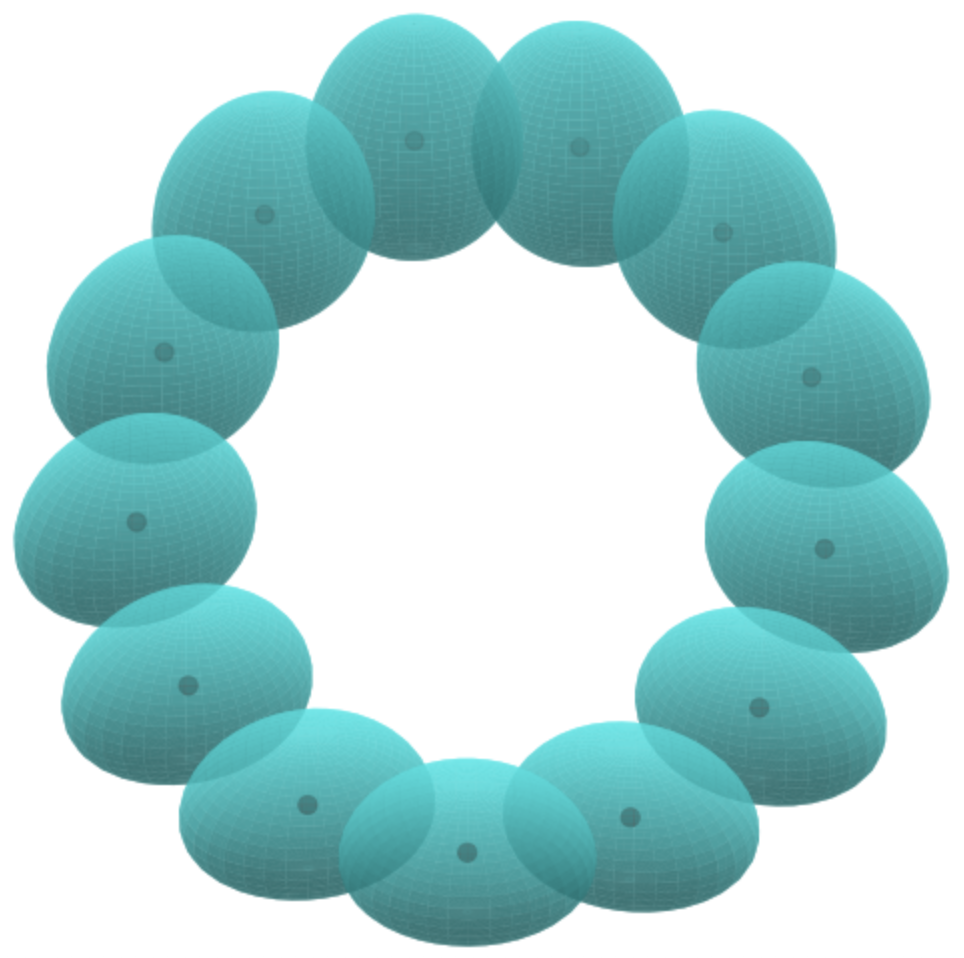}
        \caption*{ $\rho$ = 0.27}
    \end{minipage}

    \caption{The shape of the balls given by distance \eqref{distanceF}.}
    \label{fig:log-comparison2}
\end{figure}

In Figure \ref{fig:log-comparison2} we use the distance \eqref{distanceF} and the same point set lying on a circle, as in Figure \eqref{fig:log-comparison1} for illustrates the evolution of balls centered at each point for different radius. 

In the next section we will discuss more about the importance of convex and contractible sets when investigating the topology of a space whose finite point set is sampled in general position. The so-called Nerve Theorem plays a important role on this investigation and we present two different formulations for the same result in Propositions \ref{teo_nervo} and \ref{teo_nervo_2}.

%-----------------------------------------------------
\section{TDA} \label{section_tda}
%-----------------------------------------------------

\begin{defi}
    An \textbf{abstract simplicial complex} $\Delta$ is a family of subsets of a finite set $V$, such that $\Delta$ is closed under inclusion, that is, if $\sigma \in \Delta$ and $\tau\subset\sigma$ then $\tau \in \Delta$. The word \textit{abstract} can be discarded in this manuscript. $V$ is called the vertex set of $\Delta$.
\end{defi}

The elements of $\Delta$ are called faces or simplex. If an element is a subset of cardinality $k+1$, it is said to be a $k$-simplex or a face of dimension $k$. The set of faces with maximal dimension is called $facets(\Delta)$ and these elements are called \textbf{facet}. The dimension of a simplicial complex is defined to be the dimension of any facet.

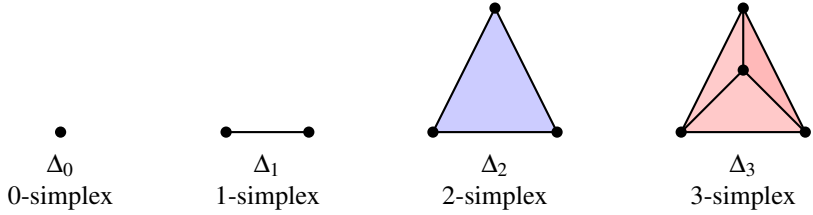
\begin{figure}[H]
\centering
\begin{tikzpicture}[scale=1.1] % Adjusted scale
% 0-simplex (a point)
\node at (-1, 0) [circle, fill=black, inner sep=1.5pt] {};
\node at (-1, -0.4) {$\Delta_0$ };
\node at (-1, -0.8) {0-simplex};
% 1-simplex (a line)
\draw[thick] (1, 0) -- (2, 0);
\node at (1, 0) [circle, fill=black, inner sep=1.5pt] {};
\node at (2, 0) [circle, fill=black, inner sep=1.5pt] {};
\node at (1.5, -0.4) {$\Delta_1$};
\node at (1.5, -0.8) {1-simplex};
% 2-simplex (a filled triangle, scaled down slightly)
\fill[blue!20] (3.5, 0) -- (5, 0) -- (4.25, 1.5) -- cycle; % Fill the triangle
\draw[thick] (3.5, 0) -- (5, 0) -- (4.25, 1.5) -- cycle; % Draw the edges
\node at (3.5, 0) [circle, fill=black, inner sep=1.5pt] {};
\node at (5, 0) [circle, fill=black, inner sep=1.5pt] {};
\node at (4.25, 1.5) [circle, fill=black, inner sep=1.5pt] {};
\node at (4.25, -0.4) {$\Delta_2$};
\node at (4.25, -0.8) {2-simplex};
% 3-simplex (a filled tetrahedron projection, scaled down slightly)
\fill[red!30, opacity=0.7] (6.5, 0) -- (8, 0) -- (7.25, 1.5) -- cycle; % Front triangle
\fill[red!20, opacity=0.7] (6.5, 0) -- (7.25, 0.75) -- (7.25, 1.5) -- cycle; % Side face
\fill[red!30, opacity=0.7] (8, 0) -- (7.25, 0.75) -- (7.25, 1.5) -- cycle; % Another face
\draw[thick] (6.5, 0) -- (8, 0) -- (7.25, 1.5) -- cycle;
\draw[thick] (6.5, 0) -- (7.25, 0.75) -- (8, 0);
\draw[thick] (7.25, 0.75) -- (7.25, 1.5);
\node at (6.5, 0) [circle, fill=black, inner sep=1.5pt] {};
\node at (8, 0) [circle, fill=black, inner sep=1.5pt] {};
\node at (7.25, 1.5) [circle, fill=black, inner sep=1.5pt] {};
\node at (7.25, 0.75) [circle, fill=black, inner sep=1.5pt] {};
\node at (7.25, -0.4) {$\Delta_3$};
\node at (7.25, -0.8) {3-simplex};
\end{tikzpicture}
\caption{\label{deaths} Examples of simplex of dimension 0,1,2 and 3. }
\end{figure}

\begin{defi}\label{nerve_defi}
    Let $\mathcal{U}=\{U_\alpha\}_{\alpha\in\Gamma}$ be a family of sets. The \textbf{nerve} of $\mathcal{U}$, denoted by $\mathcal{N(U)}$, is defined to be the simplicial complex with vertex set $\Gamma$ and $\{\alpha_1,\ldots, \alpha_{k+1}\}$ spans a $k$-simplex in $\mathcal{N(U)}$ if and only if $U_{\alpha_1} \cap \ldots \cap U_{\alpha_{k+1}} \neq \emptyset$.
\end{defi}

\begin{prop}[Nerve Theorem]\label{teo_nervo}
    Let $X$ be a topological space and $\mathcal{U}=\left\{U_\alpha\right\}_{\alpha \in \Gamma}$ a covering of $X$. If every finite intersection $U_{\alpha_1} \cap \cdots \cap U_{\alpha_k}$ is either empty or contractible, then the nerve $\mathcal{N(U)}$ is homotopy equivalent to $X$, that is,
    \[
    \mathcal{N(U)} \simeq X.
    \]
\end{prop}

\begin{remark} \label{remark2}
According to Propositions \ref{int_convex} and \ref{convex_to_contra} in the Section \ref{sec_iso_convex},  if $X\subset \mathbb{R}^n$ the nerve theorem may be modified replacing the contractibility hypothesis by convexity. That is, if $X$ is covered by closed convex sets, then the same result follows.
\end{remark}

\begin{prop}\label{teo_nervo_2}
    Let $X \subset \mathbb{R}^n$ a topological space with a closed convex covering $\mathcal{U}$. Then, $\mathcal{N(U)}$ is homotopy equivalent to $X$.
\end{prop}

The nerve theorem we presented is just one of the many versions one can find. The original discussion about nerves and nerve theorem was developed by the so-called mathematician Paul Alexandroff. A full description about nerve theorem and related results can be found in \cite{nerves_versions}. 

If $M$ is a metric space with metric $d$, we can construct a covering of $M$ by balls of ratio $\epsilon$ centered in each element $p\in M$, such that $\{\mathcal{B}_\epsilon(p)\}_{p\in M}$ is a covering of $M$. As defined just before, the nerve of this covering define a natural simplicial complex associated to the space $M$:
    
\begin{defi}
    Let $(M,d)$ be a metric space and $X$ a finite subset of $M$. Then the \textbf{Čech complex} for $X$, attached to the parameter $\epsilon$, denoted by $C(X,\epsilon)$ will be the simplicial complex whose vertex set is $X$, and where $\{p_1, \ldots, p_{k+1}\}\subset X$ spans a $k$-simplex if $\mathcal{B}_\epsilon(p_1) \cap \cdots \cap \mathcal{B}_\epsilon(p_{k+1}) \neq \emptyset $
\end{defi}

Čech complex is constructed by a high order interaction and there is no doubt that this interaction gives strong and more precise features about the topology of the point set, however its computational complexity may be a considerable setback when we look into simplexes with dimension grater than two. Looking on pairwise relation it can be note that even if $d(p_i,p_j)\leq\epsilon$, for $i,j \in \{1,\ldots,k\}$, there is no guarantee that $\mathcal{B}_\epsilon(p_1)\cap \ldots \cap\mathcal{B}_\epsilon(p_k) \neq \emptyset $, then in this case some simplexes will not be considered, i.e. the distance between a pair of points may gives less topological information than high order interaction, but on the other hand it will reduce the complexity of building simplexes of a simplicial complex way more. Just below we define the Vietoris-Rips complex of a finite set in a metric space based on pairwise relation and Proposition \ref{inclusion_1} relates both simplicial complexes.
    
\begin{defi}
    Let $(M,d)$ be a metric space and $X$ a finite subset of $M$. The \textbf{Vietoris-Rips complex} for $X$, attached to the parameter $\epsilon$, denoted by $VR(X,\epsilon)$ will be the simplicial complex whose vertex set is $X$, and where $\{p_1, \ldots, p_{k+1}\}\subset X$ spans a $k$-simplex if $d(p_i, p_j) \leq \epsilon$ for all $1 \leq i, j \leq k+1$.
\end{defi}

\begin{prop}\label{inclusion_1}
    Let $\epsilon >0$. We have the inclusions
    \[
     V R(X,  \epsilon) \subseteq \check{C}(X, \epsilon) \subseteq V R(X,  2\epsilon).
    \]
\end{prop}

\begin{proof}
    For a detailed proof see Lemma 2.13, in \cite{inclusion}.
\end{proof}

\begin{remark} 
One should note that if the vertex set $V$ of a simplicial complex $\Delta$ is contained in a metric space, and the balls $\mathcal{B}_\epsilon(p_i)$  are either closed convex or the intersection are either empty or contractible, we get the homotopy equivalence between $\Delta$ and $\mathcal{N}(\Delta)$. 
\end{remark}

Still according to the Nerve Theorem, it is intrigued that under good conditions the whole space (continuous) can be described by a simplicial complex (discrete). This bridge allows then one accesses the topology of a continuous space a simpler structure, a simplicial complex. By Definition \ref{nerve_defi} the most natural candidate to represent a continuous space is the $\check{C}ech$ complex because of the homotopy equivalence between a covering and its nerve. So then Proposition \ref{inclusion_1} provides a good approximation for the $\check{C}ech$ complex, since it sets simple bounds from lower to upper bound.

Proposition \ref{inclusion_1} goes further than being just simplicial complexes inclusions. In a more general sense, these inclusions represent a finite filtration of a set of points $X$ in $\mathbb{R}^n$. Let $K$ be a simplicial complex with vertex in $X$, we define the \textbf{filtration} of $K$ as a family of subcomplexes $(K_\alpha)_{\alpha\in \Gamma}$, such that for indices $\alpha, \beta$ in $\Gamma$, with $\alpha\leq\beta$ we have $K_\alpha\subseteq K_\beta$, and $K\subseteq\cup_{\alpha\in \Gamma}K_\alpha$. For the aim of this work we will let $\Gamma$ to be a finite set. Then, let us consider the following sequence of simplicial complexes, with their respective natural inclusions $K_r \hookrightarrow K_{r+1}$:
\[
\begin{array}{ccccccccccc}
    K_0 & \subseteq & K_1 & \subseteq & \cdots & \subseteq & K_{r-1} & \subseteq & K_r & \subseteq & \cdots \subseteq K_{n-1} \\
    \downarrow & & \downarrow & & & & \downarrow & & \downarrow &&\qquad \downarrow\\
    K_1 & \subseteq & K_2 & \subseteq & \cdots & \subseteq & K_{r} & \subseteq & K_{r+1} & \subseteq & \cdots \subseteq K_n
\end{array},
\]
where $K_n$ is the complete graph with vertexes in $X$. Given any simplicial complex $K_j$ in that sequence, one can count the number of connected component, one-dimensional holes (independent loops), two-dimensional voids or cavities and so on into higher dimensions. These quantity are the rank of abelian quotient groups called homology groups, and we denote the $k$-th homology group of $K_j$ by $H_k(K_j)$ which count the number of independent $k$-dimensional holes in $K_j$. Due to homology be functorial, each inclusion $i_r : K_r \hookrightarrow K_{r+1}$ induces a homology map $l_{i_r}^k:H_k(K_r)\rightarrow H_k(K_{r+1})$ \cite{hatcher}, so this sequence of homology groups forms a diagram:
\[
H_k\left(K_0\right) \longrightarrow H_k\left(K_1\right) \longrightarrow \cdots \longrightarrow H_k\left(K_n\right).
\]

This diagram allows to track topological features such as independent loops or voids along the filtration $K_j\subseteq K_{j+1}$. For instance, if $\beta$ is $k$-dimensional independent hole which arises in $H_k(K_i)$, survives into $H_k(K_i)$ but disappear in $H_k(K_{i+1})$, we say that $\beta$ had $[i,j)\subset\Gamma$ time of life, birth in $i$ and death in $j$, that is the persistence of $\beta$. Hence, the homology map diagram allows us to measure birth and death of holes, looking at the image of each map $l_{i_r}^k$ and the possible composition, lets say $l_{i_s}^k\circ\cdots\circ l_{i_r}^k$, if a hole appears in $H_k(K_r)$ and disappears in $H_k(K_s)$. The pioneering work on tracking how long a hole persist in the filtration is \cite{tda2}, in which the authors call it \textbf{persistent homology}. For each feature in the filtration, whose persist from $b$ to $d$, the point $(b, d) \in \mathbb{R}^2$ corresponds to the birth and death of this individual homology classes in $H_k$. Hence, one obtains a persistence diagram for each dimension $k$, whose together obtain the so-called \textbf{persistent diagram} of the filtration. 

In the paper \cite{nobrega1} we can view a filtration process on finite set of point, where the threshold $\epsilon$ increases, and then edges, triangles and tetrahedrons are built, and on the top of the simplicial complexes are the sublevel sets of a torus to represent the same idea of birth and death topological features. Through this process it is possible to see that at some stages, appears an annulus and die a one-independent loop, making this an important mark in the filtration. At some times the persistent diagram describes the birth and death of classes in the homology groups of dimensions $0$, $1$ and $2$. 

%-----------------------------------------------------
\section{Finsler Structures in TDA}\label{section_fisnler_tda}
%-----------------------------------------------------

TDA fundamentally relies on classical tools to investigate the topology of a continuous space $\Omega$, as well as the topology of a discrete set $X \subset \Omega$, through simplicial complexes such as the Vietoris–Rips complex and the Čech complex, which are built using distance functions in the context of a metric space. Based on this, it is possible to infer a topological structure from a point cloud and these approaches are well established and have proven effective in many applications. In this section, we present a framework that has the potential to be able to capture richer structural features beyond what standard metrics can reveal.

Let $(M,F)$ be a Finsler space, by Remark \ref{remark1} for each $p \in M$, the function $F(p,v)$ defines a quasi-norm for $v$ in $T_pM$. For simplicity, once set the point $p$, lets write $F_p(v)$ instead of $F(p,v)$ to make clear that $F$ is acting only on the vectors $v\in T_pM$.

Since $F_p:T_pM \rightarrow \mathbb{R}$ is a quasi-norm, $F_p$ induces then two semi-metric, or asymmetric distances, on the product $T_p M \times T_p M$ as follow:
\begin{align*}
    d^{+}_p: T_p M \times T_p M &\rightarrow \mathbb{R}, \quad (u,v) \mapsto d^{+}_p(u,v) = F_p(v - u) \\
    d^{-}_p: T_p M \times T_p M &\rightarrow \mathbb{R}, \quad (u,v) \mapsto d^{-}_p(u,v) = F_p(u - v).
\end{align*}

Let $\epsilon>0$, and fix $v\in T_pM$, $F_p$ also generates two balls of radius $\epsilon$ centered at $v$. These balls are namely the forward $\mathcal{B^{+}}(v,\epsilon)$ and backward $\mathcal{B^{-}}(v,\epsilon)$ balls:
\begin{align*}
    \mathcal{B^{+}}(v,\epsilon) = \{w\in T_pM ; d^{+}_p(v,w) \leq \epsilon  \} \\
    \mathcal{B^{-}}(v,\epsilon) = \{w\in T_pM ; d^{-}_p(v,w)) \leq \epsilon  \}.
\end{align*}

It is very simple to see that if $F$ is absolute homogeneous, that is $F_p(v)=F_p(-v)$, then $d_p^{+}=d_p^{-}$ and consequently $\mathcal{B}^{+}(\nu, \boldsymbol{\varepsilon})=\mathcal{B}^{-}(\nu, \boldsymbol{\varepsilon})$. Letting $p$ varies on $M$ we obtain two families of semi-metrics and their respective family of balls.

For now, let us restrict the space $M$ to a convenient set which will be way more useful for our purposes. Let $\Omega$ be an open connected subset of $\mathbb{R}^n$ consider $M=\Omega$ and $\mathbb{F}=(F,\Omega$) a Finsler space. In this case, the tangent space of $\Omega$ can be seen as $T_pM \cong  \mathbb{R}^n$ and the tangent bundle as a Cartesian product $T \Omega \cong \Omega \times \mathbb{R}^n$. Therefore, two points $p,q$ in $\Omega$ will be seen as vectors of $T_pM$ and we can calculate the distance between $p$ and $q$ as $d_p(p,q) = F_p(q-p)$. This makes it possible to build forward and backward balls centered in each point $p$ of $\Omega$ with the corresponding distance associated to $p$. 

Now we have two types of Finsler distances that can be considered in $\Omega$: the one defined in this section, given by the Finsler norm $F_p$, and the one defined in Section \ref{section_finsler}, obtained by taking the infimum of all curves connecting two points. In fact, under good conditions, they can be equivalent; this is shown in the next result.

\begin{prop}\label{pointdistance}
    If $\Omega \subset \mathbb{R}^n$ and $F$ is a Finsler metric that does not depends of first variable then 
    \[
    d_F(x,y) = F(x,y-x).
    \]
\end{prop}
\begin{proof}
    If $F(p,v)$ does not depend on first variable then the curves that perform the distance $d_F$ satisfies $\displaystyle\dfrac{d^2 p^i}{ds^2}=0,$ i.e., the curves are $\gamma^i(s) = a_i~s + b_i$ for $a_i, b_i$ constants. If $\gamma^i(0) = x^i$, $\gamma^i(1)=y^i$, we have $a_i = y^i - x^i$ and $b^i = x^i$, then
    \begin{align}
    d_F(x,y) = \int_0^1 F(a_i~s + b_i, a_i)~ds = \int_0^1 F(s(y^i - x^i) + x^i, y^i - x^i)~ds \\
 = \int_0^1 F(y^i - x^i)~ds= s~F(y^i - x^i) \Big|_0^1 = F(y-x) = F_x(y-x) .
     \end{align}
\end{proof}

\begin{remark}
From here, we will work only on the forward balls but the same ideas hold for the backward balls, and for simplicity we will avoid the symbol $"+"$ on the top of $\mathcal{B}^{+}$ and $d_p^{+}$. 
\end{remark}

Since we can identify any point $p$ of $\Omega$ as a vector in $\mathbb{R}^n$, we obtain a covering for $\Omega$ as the family $\mathcal{B}=\{\mathcal{B}(p, \boldsymbol{\varepsilon})\}_{p\in \Omega}$. For such distances to be used in TDA, it is important that the generated balls satisfy the assumptions of the Nerve Theorem, we then have the following result:

\begin{prop}
    Each ball in $\mathcal{B}$ is convex, and hence all non-empty finite intersection of balls in $\mathcal{B}$ is contractible.
\end{prop}
\begin{proof}
    As one of the main properties of $F$ is the positive definiteness of the tensor $g_{i j}(p, v)=\frac{1}{2} \frac{\partial^2 F^2}{\partial v^i \partial \nu^j}$, then this implies that $F_p(v)$ satisfies the triangle inequality. This is sufficient to obtain the convexity on the balls of $\mathcal{B}$. Each intersection of balls is either empty or convex, then the result follows by Proposition \ref{convex_to_contra}.
\end{proof} 

\begin{figure}[H]
    \centering
    \begin{minipage}[t]{0.19\textwidth}
        \centering
        \includegraphics[width=\linewidth]{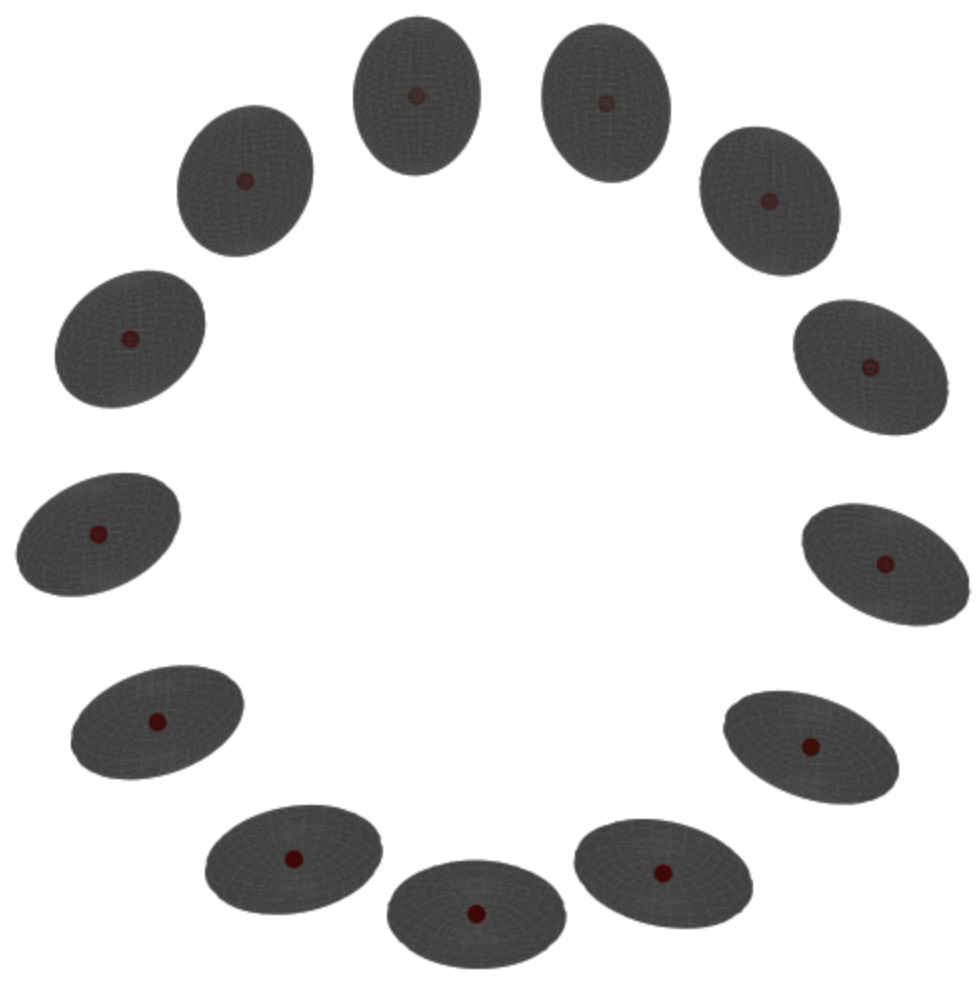}
        \caption*{$\rho$ = 0.08}
    \end{minipage}
    \hfill
    \begin{minipage}[t]{0.19\textwidth}
        \centering
        \includegraphics[width=\linewidth]{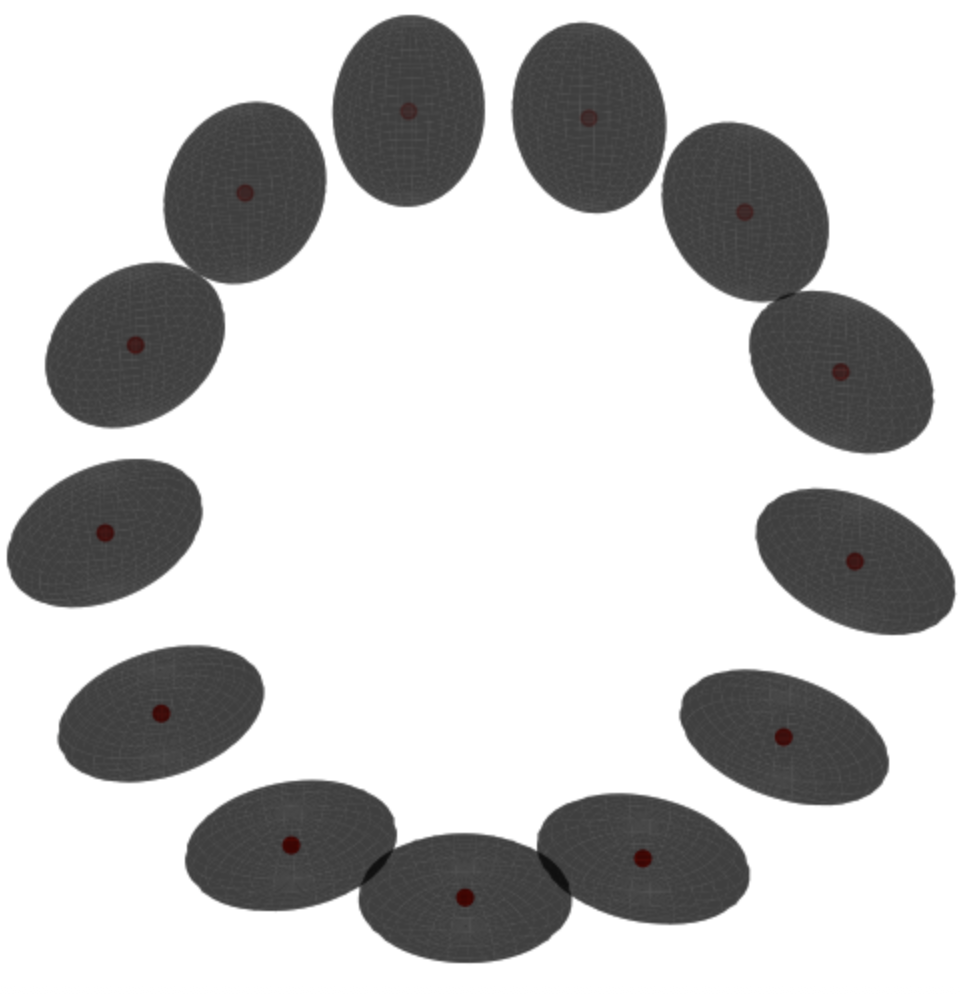}
        \caption*{ $\rho$ = 0.1 }
    \end{minipage}
    \hfill
    \begin{minipage}[t]{0.19\textwidth}
        \centering
        \includegraphics[width=\linewidth]{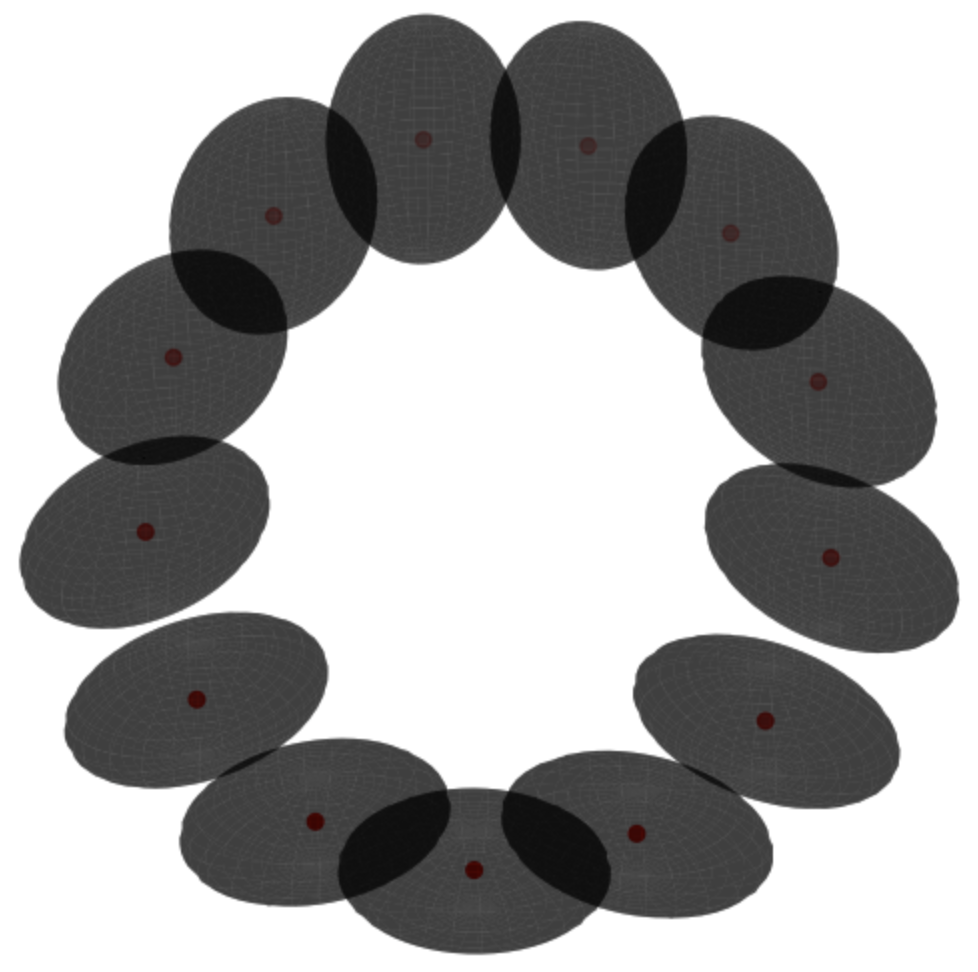}
        \caption*{ $\rho$ = 0.14}
    \end{minipage}
    \hfill
    \begin{minipage}[t]{0.19\textwidth}
        \centering
        \includegraphics[width=\linewidth]{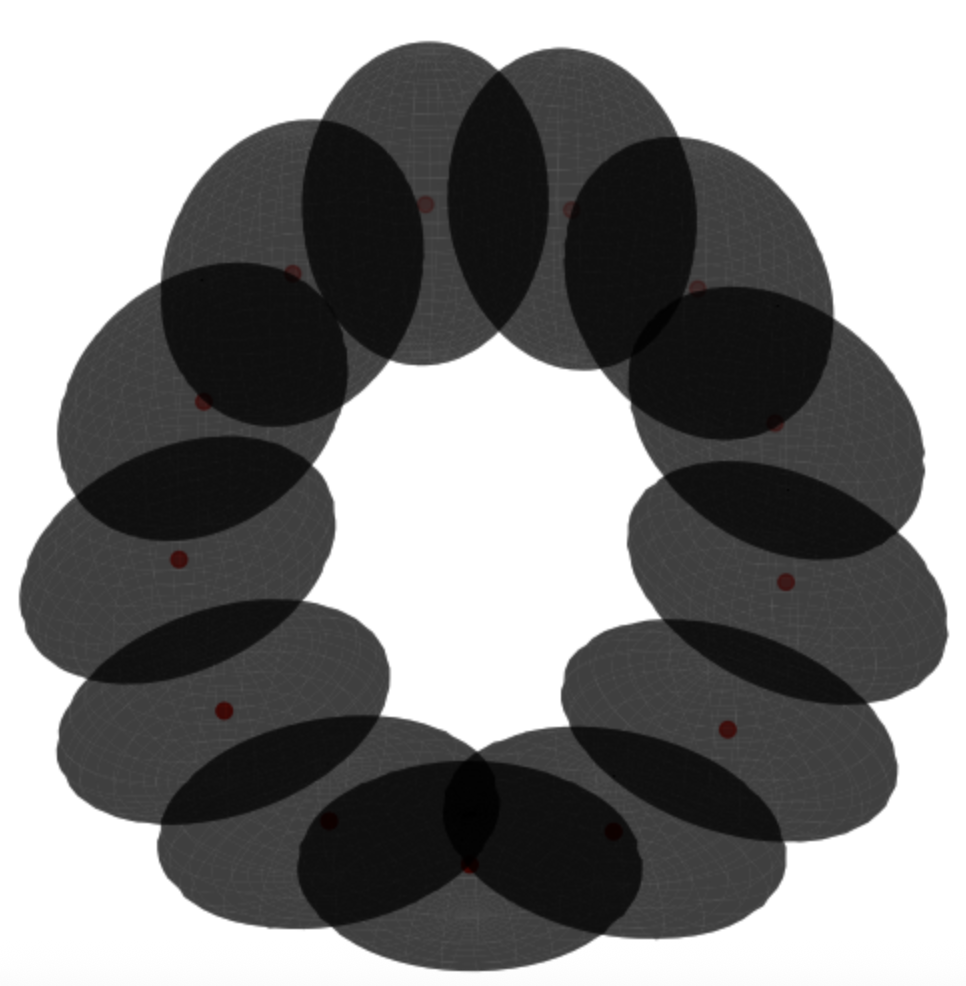}
        \caption*{ $\rho$ = 0.2}
    \end{minipage}

    \caption{The shape of the balls given by Finsler metric $F_p(v)=\left(\sum_{i=1}^n\left(\frac{v^i}{p^i}\right)^2\right)^{\frac{1}{2}}$.}
    \label{fig:log-comparison3}
\end{figure}

In Figure \ref{fig:log-comparison3} we again use the set of points lying on a circle, as in Figures \eqref{fig:log-comparison1} and \eqref{fig:log-comparison2} to illustrate the evolution of forward balls centered at each point for different radius, but now we are using the Finsler metric \eqref{metrica_raiz} with the method just introduced above. Note that in this case, the metric $F_p$ is absolute homogeneous and the forward balls and the backward balls are the same.

Let $X\subset\Omega$ be a finite point set in general position. There is a natural way to generalize the Vietoris-Rips complex in this Finsler space. The \textbf{generalized Vietoris-Rips complex} in $\mathbb{F}$ associated to $X$ and $\epsilon>0$, denoted by $Rips(X,\epsilon)$, is the simplicial complex such that $\{p_1, \ldots, p_{k+1}\}$ spans a $k$-simplex if
\[
d_{p_i}(p_i, p_j)\leq \epsilon \quad \text{for all} \quad p_i, p_j \in \left\{p_1, \ldots, p_{k+1}\right\}.
\]
On the other hand the \textbf{generalized $\check{C}$ech complex} associated to $X$ and $\epsilon>0$, denoted by $\check{C}ech(X,\epsilon)$, is the simplicial complex such that $\left\{p_1, \ldots, p_{k+1}\right\}$ spans a $k$-simplex if
\[
\bigcap\mathcal{B}(p_i,\epsilon)\neq \emptyset \quad \text{for all} \quad p_i\in \left\{p_1, \ldots, p_{k+1}\right\}.
\]

Alternatively, $\check{C}ech(X,\epsilon)$ can be thought as the nerve $\mathcal{N(B)}$. Lets say that $\mathcal{B}_\epsilon$ is the family of all the balls of radius $\epsilon$ centered at $p\in X$ as above and $\mathcal{B}_\delta$, with $\delta>\epsilon$, the respective family following the same thoughts. If $\delta-\epsilon$ is sufficiently small, then there are no new intersections over the balls, and hence $\mathcal{N(B_\epsilon)}=\mathcal{N(B_\delta)}$. So that $\check{C}ech(X, \epsilon) = \check{C}ech(X, \delta)$.

\begin{remark}\label{remark7}
One can easily check that $Rips(X,\epsilon) \subseteq \check{C}ech(X, \epsilon)$ still holds. If we have that $\left\{p_1, \ldots, p_{k+1}\right\}$ spans a $k$-simplex of $Rips(X,\epsilon)$, then $d_{p_i}\left(p_i, p_j\right) \leq \varepsilon$, for all $p_i, p_j \in \{p_1, \ldots, p_{k+1}\}$. Choose $p$ as $p_1$, for instance, then $d_{p_i}\left(p_i, p\right) \leq \varepsilon$, for all $p_i \in \{p_1, \ldots, p_{k+1}\}$. It means that $p\in \bigcap \mathcal{B}\left(p_i, \varepsilon\right)$.
\end{remark}

Ir order to justify the use of Finsler spaces in TDA, given a finite point set in general position we present a new method of building simplicial complex based on pairwise relation between the points, but it naturally carries more information about the points beyond just the distance between them, so that it becomes richer in information than the generalized Vietoris-Rips complex. 

\begin{defi}\label{new_complex}
    Let $X$ be a finite point set in a Finsler space $(\Omega,F)$. We denote by $\Delta(X,\epsilon)$ the simplicial complex associated to $X$ and $\epsilon>0$ such that $\{p_1,\ldots,p_{k+1}\}$ spans a $k$-simplex of $\Delta(X,\epsilon)$ if:
    \[
    d_{p_i}\left(p_i, p_j\right) + d_{p_j}\left(p_j, p_i\right) \leq 2\varepsilon \quad \text{for all} \quad p_i, p_j \in \left\{p_1, \ldots, p_{k+1}\right\}.
    \]
\end{defi}

\begin{prop}\label{inclusion_2}
    Let $\epsilon>0$ and $X$ be a finite point set in a Finsler space $(\Omega,F)$. Then,
    \[
    \text{Rips}(X,\epsilon)\subseteq\Delta(X,\epsilon)\subseteq \text{Rips}(X,2\epsilon).
    \]
\end{prop}

\begin{proof}
For the first inclusion, it is sufficient to note that if the subset $\left\{p_1, \ldots, p_{k+1}\right\}$ spans a $k$-simplex in $\text{Rips}(X,\epsilon)$, then $d_{p_i}\left(p_i, p_j\right) \leq \epsilon$, for all $p_i, p_j \in \{p_1, \ldots, p_{k+1}\}$. Since the this condition holds for any $p_i$ and $p_j$ we can exchange $i$ with $j$ so that we also have $d_{p_j}\left(p_j, p_i\right) \leq \epsilon$. Hence, by the positivity of each $d_{p_i}$ we have the inequalities
    \[
    d_{p_i}\left(p_i, p_j\right) \leq d_{p_i}\left(p_i, p_j\right)+d_{p_j}\left(p_j, p_i\right) \leq 2 \varepsilon, \quad \text { for all } \quad p_i, p_j \in \{p_1, \ldots, p_{k+1}\}
    \]
Conversely, using the definitions of $\text{Rips}(X,\epsilon)$ and $\Delta(X,\epsilon)$, and the same argument on the positivity of $d_{p_i}$, we can prove that $d_{p_i}\left(p_i, p_j\right) \leq 2\epsilon$ for all $p_i, p_j \in \{p_1, \ldots, p_{k+1}\}$.
\end{proof}

An important result in TDA, when you consider the point set into a space with one symmetric distance, is the inclusion given by Proposition \ref{inclusion_1} which approximates the $\text{Č}$ech complex by Vietoris-Rips complex. Proposition \ref{inclusion_2} establishes the same kind of inclusions with lower and upper bounds given by Vietoris-Rips complexes. In the following result we guarantee that \text{Č}ech complexes are included in Vietoris-Rips complexes of finite radius, but not necessarily less than or equal to $2\epsilon$.

\begin{prop}\label{inclusion_3}
    Let $\epsilon>0$ and $X$ be a finite point set in a Finsler space $(\Omega,F)$. Then,
    \[
    \text{Rips}(X,\epsilon)\subseteq\check{C}ech(X,\epsilon)\subseteq \text{Rips}(X,\epsilon +\delta),\quad \delta<\infty.
    \]
\end{prop}
\begin{proof}
    The first inclusion was mentioned in the Remark \ref{remark7}. For the other inclusion, let $\left\{p_1, \ldots, p_{k+1}\right\}$ spans a $k$-simplex in $\check{C}ech(X, \varepsilon)$. By definition of $\check{C}ech(X,\epsilon)$, we can choose $y\in \bigcap \mathcal{B}\left(p_i, \varepsilon\right)$, with the property $d_{p_i}\left(p_i, y\right)\leq \epsilon$, for all $p_i \in \{p_1, \ldots, p_{k+1}\}$. By triangle inequality 
    \[
    d_{p_i}\left(p_i, p_j\right)\leq d_{p_i}\left(p_i, y\right) +d_{p_i}\left(y, p_j\right)\leq \varepsilon +\delta,
    \]
    where $\delta = \displaystyle\sup_{p_i,p_j}d_{p_i}(y,p_j)$. Since all norms in $\mathbb{R}^n$ are equivalents, and then they can be comparable, and beside that $d_{p_i}\left(y, p_j\right)=F_{p_i}(y-p_j)$, $\delta$ must be a real number.
\end{proof}

\begin{remark}
Now, we will present a type of Finsler metric that generates an asymmetric distance, which can be used in the theory we are developing.
\end{remark} 

% \begin{defi}
% A Finsler metric $F(p,v) = \alpha(p,v) + \beta(p,v)$ is called a \textbf{Randers metric} when $\alpha(p,v)$ is a Riemannian metric $\sqrt{a_{ij}(p)v^iv^j}$ in $M$ and $\beta(p,v)$ is a differential 1-form $b_i(p)v^i$ in $M$. The Finsler space $F = (M, F = \alpha + \beta)$ with a Randers metric is called a \textbf{Randers space}.
% \end{defi}

Let $f: \Omega\subseteq\mathbb{R}^n \rightarrow \mathbb{R}$ a non-negative smooth function and $\omega$ a 1-form in $\Omega$ with coefficients $b_i(p)=\omega_p \left(\frac{\partial}{\partial p^i}\right) \in \mathbb{R}$. Define the function $F(p,v)$ by 
\begin{equation}\label{alphametric}
    F(p,v) = \left( \sum_{i=1}^n f_i(p) | v^i |^\alpha \right)^\frac{1}{\alpha} + \sum_{i=1}^n b_i(p)\big(f_i(p)\big)^\frac{1}{\alpha}v^i. 
\end{equation}

The first term of \eqref{alphametric} is a weighted $\alpha$-norm. One should note that this norm by itself defines a Finsler metric $\alpha_f(p,v)=\left(\sum_{i=1}^n f_i(p)\left|v^i\right|^\alpha\right)^{\frac{1}{\alpha}}$. Before we prove this function is indeed a Finsler metric, lets see its balls intersection behavior as we did in Figures \eqref{fig:log-comparison1}, \eqref{fig:log-comparison2} and \eqref{fig:log-comparison3}. Here we use the same points lying on a circle:
\begin{figure}[H]
    \centering
    \begin{minipage}[t]{0.19\textwidth}
        \centering
        \includegraphics[width=\linewidth]{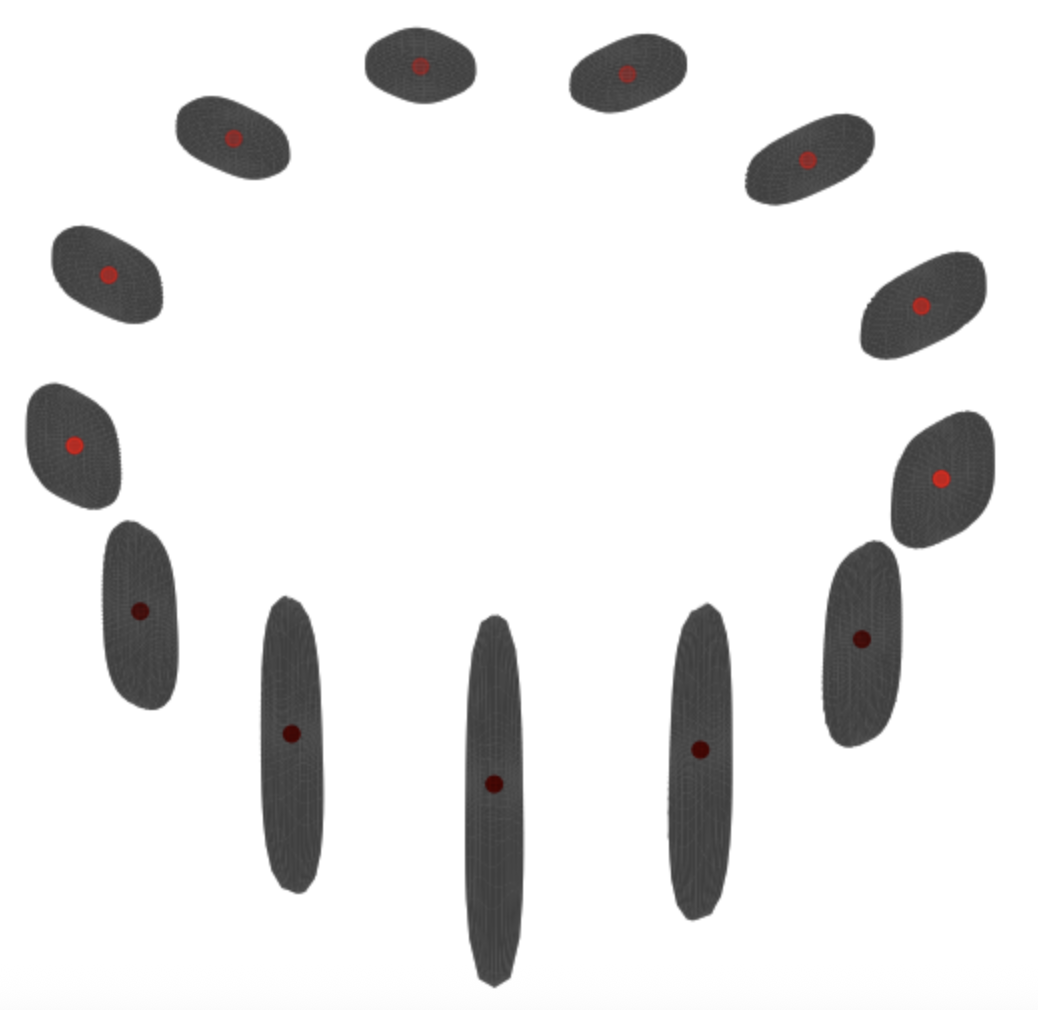}
        \caption*{$\rho$ = 100}
    \end{minipage}
    \hfill
    \begin{minipage}[t]{0.19\textwidth}
        \centering
        \includegraphics[width=\linewidth]{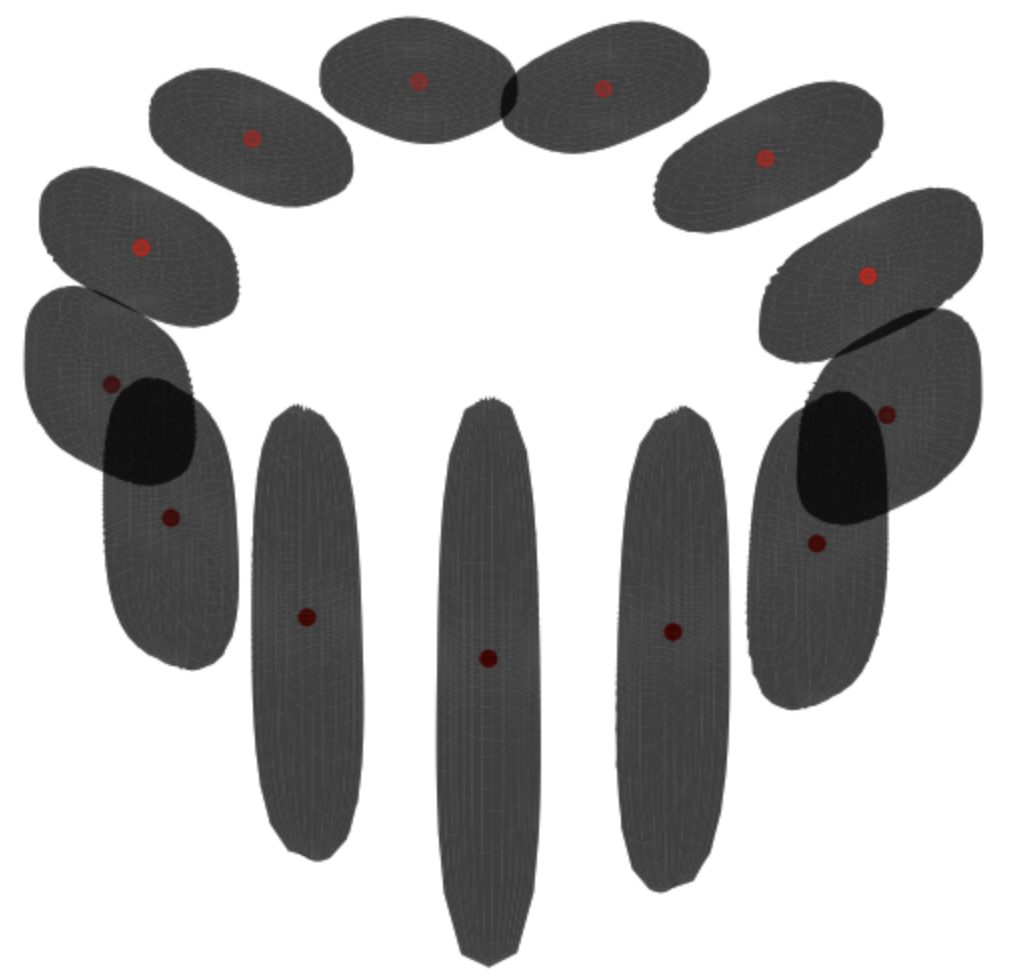}
        \caption*{ \small$\rho$ = 200 }
    \end{minipage}
    \hfill
    \begin{minipage}[t]{0.19\textwidth}
        \centering
        \includegraphics[width=\linewidth]{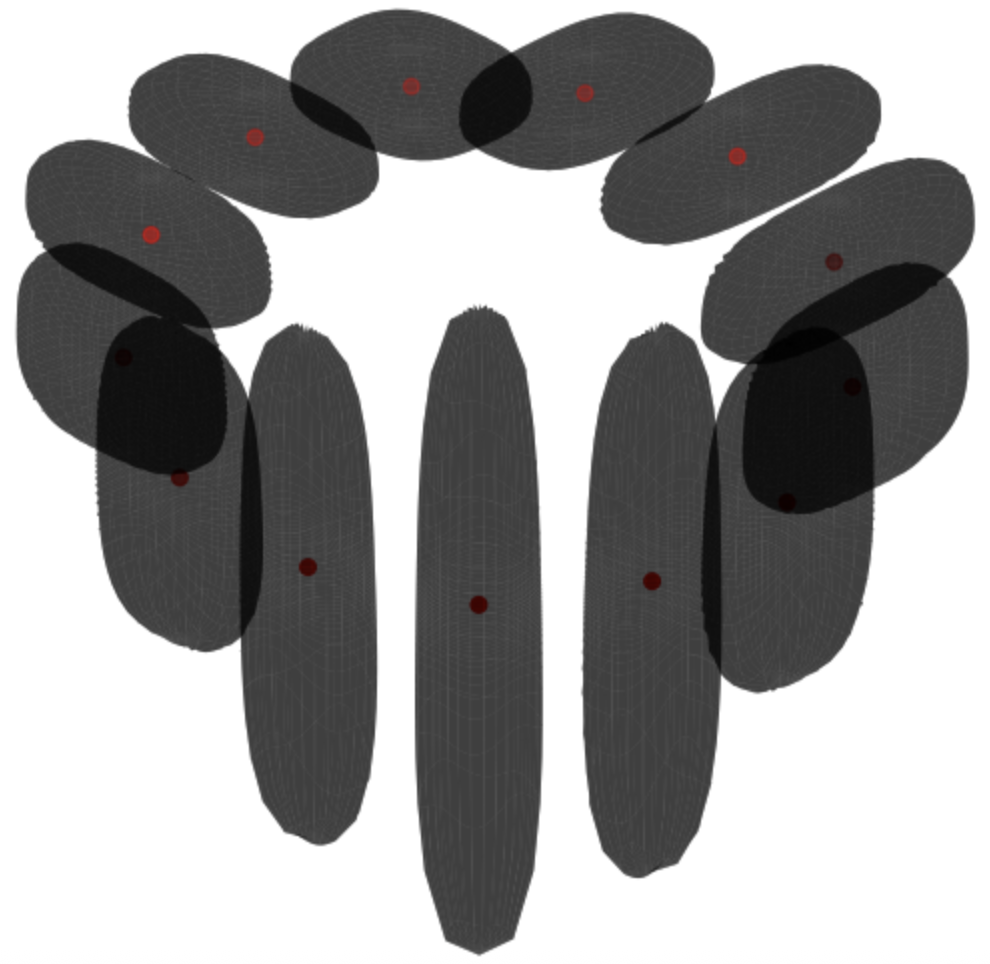}
        \caption*{ \small$\rho$ = 260}
    \end{minipage}
    \hfill
    \begin{minipage}[t]{0.19\textwidth}
        \centering
        \includegraphics[width=\linewidth]{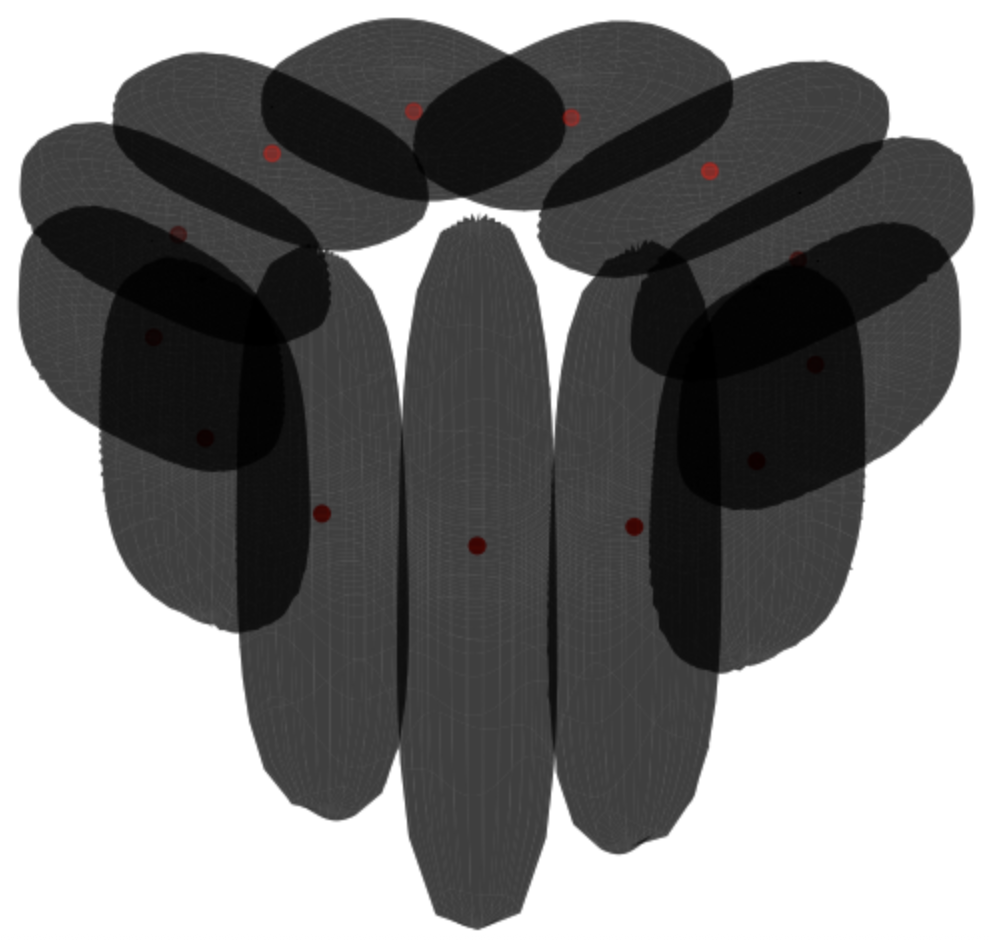}
        \caption*{ \small$\rho$ = 360}
    \end{minipage}

    \caption{The shape of the balls given by Finsler metric $F_p(v)=\left( \sum_{i=1}^n p_{i}^2 | v^i |^3 \right)^\frac{1}{3} + \sum_{i=1}^n b_i(p)p_{i}^{\frac{2}{3}}v^i$, with $b_i(p) = \frac{p_i}{1+\|p\|_{3/2}}$.}
    \label{fig:log-comparison4}
\end{figure}
\begin{prop}\label{our_metric}
    If $\alpha$ and $\beta$ are conjugate exponents, $1<\alpha<\infty$ and $\omega$ is a 1-form such that the vector $b=(b_i(p)) \in T_p\Omega$ satisfies $\|b\|_\beta < 1$, then $F$ given by \eqref{alphametric} is a Finsler metric.
\end{prop}

\begin{proof}
    It is clear that $F$ is positively homogeneous of degree 1 in $v$ and for $1<\alpha<\infty$, the function $F$ is differentiable out of the null-section of $T\mathbb{R}^n$. 
    
    Consider $\alpha_f(p,v) = \displaystyle \left( \sum_{i=1}^n f_i(p) | v^i |^\alpha \right)^\frac{1}{\alpha}$. So that, for each $p\in M$, $\alpha_f(p,v)$ defines a $\alpha$-norm of the vector $u = \big(f_i(p)^\frac{1}{\alpha} v^i\big) \in T_pM$. The Cauchy-Schwarz inequality guarantee 
    \[
    |b_i(p)\big(f_i(p)\big)^\frac{1}{\alpha}v^i| \leq \|b\|_\beta \|u\|_\alpha,
    \]

Now, note that $F$ provides the inequality
    \begin{equation*}
    \begin{aligned}
        F(p,v) &\geq \alpha_f(p,v) - \|b\|_\beta \|u\|_\alpha\\
        &= \alpha_f(p,v) - \|b\|_\beta \left(\sum_{i=1}^n f_i(p)|v^i|^\alpha\right)^\frac{1}{\alpha}\\
        &= \alpha_f(p,v) \big( 1 - \|b\|_\beta\big). 
    \end{aligned}
    \end{equation*}
    
    Therefore, $F$ is positive definite for $\|b\|_\beta < 1$. Since $F$ is strictly convex, follows that $g_{i j}(p, v)=\displaystyle\frac{1}{2} \frac{\partial^2 F^2}{\partial v^i \partial v^j}$ is positive definite.
\end{proof}

In general, given an asymmetric distance it is possible to define two $\check{C}$ech complexes, one given by the nerve of forward balls and the other built as the nerve of backward balls. The same holds for $Rips$ complexes and the simplicial complex defined in Definition \ref{new_complex}. Let $\check{C}ech^+(X,\epsilon)$ and $\check{C}ech^-(X,\epsilon)$ be the $\check{C}$ech complexes associated to forward and backward balls of $d_p^F(x, y)$ respectively. According to the Finsler metric proposed in Proposition \ref{our_metric}, the following result establishes good relation between the $\check{C}$ech given by the forward and backward balls

\begin{prop}
For the Finsler metrics given by (\ref{alphametric}), we have
\[
\check{C}ech^-(X, \varepsilon) \subseteq \check{C}ech^+\left(X,\frac{\varepsilon (1+c)}{1-c}\right)\subseteq \check{C}ech^-\left(X,\frac{\varepsilon (1+c)^2}{(1-c)^2}\right),\quad \|b\|_\beta \leq c<1.
\]
\end{prop}
\begin{proof}
As proved in Proposition $\ref{our_metric}$, we already have $(1-\|b\|_\beta) \alpha_f(p, v) \leq F(p, v)$. But now, note that $F(p,v)\leq (1+\|b\|_\beta)\alpha_f(p,v)$ because of the triangle inequality applied on the absolute value of $F(p,v)$. Taking $c$, $\|b\|_\beta \leq c<1$, then
\begin{equation*}\label{alphainequality1}
    (1-c)\alpha_f(p, v)\leq F(p,v)\leq (1+c) \alpha_f(p, v).
\end{equation*}
So then, for each $p$, the inequality above induces a comparison between the point-wise distances $d^F_p$ and $d^{\alpha_f}_p$ for $F$ and $\alpha_f$ respectively.
\begin{equation*}\label{alphainequality2}
(1-c) d^{\alpha_f}_p(x,y) \leq d^F_p(x,y) \leq(1+c) d^{\alpha_f}_p(x,y), \quad \text{for all}, \quad x,y\in\Omega.
\end{equation*}
So if  $\sigma \in \check{C}ech^-(X, \varepsilon)$, then the first inequality guarantee $\sigma \in \check{C}ech_{\alpha_f}\left(X,\frac{\varepsilon}{1-c}\right)$ and by the second inequality we also have $\sigma \in \check{C}ech^+\left(X,\frac{\varepsilon (1+c)}{1-c}\right)$. The first inclusion is then obtained: 
\[\check{C}ech^-(X, \varepsilon) \subseteq \check{C}ech^+\left(X,\frac{\varepsilon (1+c)}{1-c}\right).\]
In the same way, if $\sigma \in \check{C}ech^+\left(X, \frac{\varepsilon (1+c)}{1-c}\right)$ then $\sigma \in \check{C}ech_{\alpha_f}\left(X,\frac{\varepsilon (1+c)}{(1-c)^2}\right)$ and this implies $\sigma \in \check{C}ech^-\left(X,\frac{\varepsilon (1+c)^2}{(1-c)^2}\right)$, which gives the second inclusion: 
\[\check{C}ech^+\left(X, \frac{\varepsilon (1+c)}{1-c}\right) \subseteq \check{C}ech^-\left(X,\frac{\varepsilon (1+c)^2}{(1-c)^2}\right).\]
\end{proof}

%-----------------------------------------------------
\section{Stability of the persistent diagrams}\label{section_stability}
%-----------------------------------------------------

Motivated by the paper \cite{chazal_stability} we also want to give a good satisfaction about the persistence diagram associated to the filtration of two simplicial complexes with vertex sets $X$ and $Y$, when we use a Finsler metric to define the simplexes of each simplicial complex. So throughout this section we make use of definitions and results of theorems extensively explained in \cite{chazal_defis, chazal_results}. 

Let $(K_\epsilon)_{\epsilon\in \mathbb{R}}$ the filtration of a fixed simplicial complex $S$, the family of vector spaces $(H(K_\epsilon))_{\epsilon\in \mathbb{R}}$ is a \textbf{persistence module} over $\mathbb{R}$. Let $X,Y$ be finite subsets of $\mathbb{R}^n$, the $\text{Č}$ech complexs and Vietoris-Rips complexs of $X$ and $Y$ are family of simplicial complexs, and your homology groups are persistence modules equipped with linear maps given by inclusions that have finite rank. These persistence modules associated to $X$ and $Y$ are $\delta$-interleaved when the correspondence $C: X \rightrightarrows Y$ and $C^T: Y \rightrightarrows X$ (see \cite{chazal_stability}) are $\delta$-simplicial, this happens when given two filtered simplicial complex $\mathbb{S}$ and $\mathbb{T}$ with vertex set $X$ and $Y$, for any $\epsilon \in \mathbb{R}$ and any simplex $\sigma \in S_\epsilon$, every finite set of $C(\sigma)$ is a simplex of $T_{\epsilon + \delta}$ and the same happens of $C^T$. 

The \textbf{Gromov-Hausdorff distance} for $X$ and $Y$ can be write as
\[
\mathrm{d}_{\mathrm{GH}}(X, Y)=\frac{1}{2} \inf \{\operatorname{dis}(C): C \text { is a correspondence } X \rightrightarrows Y\},
\]
where
\[
\operatorname{dis}(C)=\sup \left\{\left|d_X\left(x, x^{\prime}\right)-d_Y\left(y, y^{\prime}\right)\right|:(x, y),(x^{\prime}, y^{\prime}) \in C\right\}.
\]

Let $U$ and $V$ be two persistence modules associated with filtration of $X$ and $Y$, respectively. The \textbf{bottleneck distance} between the persistence diagrams $dgm(U)$ and $dgm(V)$ is 
\[
d_b (dgm(U), dgm(V)) = \inf_\eta  ~\sup_x  \parallel x - \eta(x) \parallel_\infty , 
\]
where the supremum is to take over every $x \in dgm(U)$ and the infimum is to take over every bijection $\eta$ between $dgm(U)$ and $dgm(V)$. 

The main result about stability of diagrams is the following theorem.

\begin{theo}
    Let $X$, $Y$ be bounded metric spaces. Then
\[
\begin{aligned}
\text{d}_b(\text{dgm}(H(VR(X))),\text{dgm}(H(VR(Y))))\leq \mathrm{d}_{\mathrm{GH}}(X, Y) \\
\text{d}_b(\text{dgm}(H(\check{C}(X))),\text{dgm}(H(\check{C}(Y))))\leq \mathrm{d}_{\mathrm{GH}}(X, Y)
\end{aligned}
\]
\end{theo}
\begin{proof}
The full and detailed proof of this result is contained in \cite{chazal_stability}.
\end{proof}

This result plays a suitable hole on understanding how topological features of a simplicial complex differ from another when looking into their vertex sets. One may note that from the right side, $\mathrm{d}_{\mathrm{GH}}$ controls how $X$ differs from $Y$. Since the left side above controls the difference between topological features of each filtration, this theorem says that $\text{Rips}$ and $\check{C}ech$ complexes do not modify too much their topological invariants under small perturbations.

The approach in Section \ref{section_fisnler_tda} deals with a new, even the same, method of building simplicial complexes using a more general notion of distance, namely given by a Finslerian metric. One may ask if that filtration process guarantee also stability on the topological features for small perturbation on the points. Yes, the answer is affirmative. The closer the sets \( X \) and \( Y \) are in the Gromov--Hausdorff distance, the smaller the difference in their topological features captured during the filtration process. In particular, the bottleneck distance between the persistence diagrams of \( X \) and \( Y \) is bounded above by \( d_{\mathrm{GH}}(X, Y) \), reflecting the stability of persistent homology under perturbations in the underlying Finsler space.

\begin{theo}
Let $X$, $Y$ be bounded sets of $\mathbb{R}^n$ and $F$ a Finsler metric defined on an open set contaning $X$ and $Y$. Then,
\[
\begin{aligned}
\text{d}_b(\text{dgm}(H(\text{Rips}(X))),\text{dgm}(H(\text{Rips}(Y))))\leq \mathrm{d}_{\mathrm{GH}}(X, Y) \\
\text{d}_b(\text{dgm}(H(\check{C}ech(X))),\text{dgm}(H(\check{C}ech(Y))))\leq \mathrm{d}_{\mathrm{GH}}(X, Y) \\
\text{d}_b(\text{dgm}(H(\Delta(X))),\text{dgm}(H(\Delta(Y))))\leq \mathrm{d}_{\mathrm{GH}}(X, Y)
\end{aligned}
\]
\end{theo}
\begin{proof}
The result follow easily, once made an adaptation in the definition of distortion of a correspondence. For us, not getting into the details, the distortion of $C$ under $F$ is the set 
\[
\operatorname{dis}(C)=\sup \left\{\left|d_x\left(x, x^{\prime}\right)-d_y\left(y, y^{\prime}\right)\right|:(x, y),\left(x^{\prime}, y^{\prime}\right) \in C\right\}.
\]
Note that the distance in $X$ and $Y$ is induced as the same way as in Section \ref{section_fisnler_tda}.
\end{proof}

\section{Discussion and Conclusion}

\par Our discussion in Section \ref{sec_iso_convex}  suggest that distances derived from classical information theoretic divergences can be replaced by geodesic distances induced by suitable Finsler metrics. In this perspective, data points originally treated as elements of $\mathbb{R}^n$ may instead be viewed as lying on a differentiable surface (or more generally, a Finsler manifold) endowed with a geometry that reflects the intrinsic structure of the dataset. We established a method for building simplicial complexes with a Finsler metric that
incorporates both position and direction which always produces convex 
balls in a possible filtration. Consequently, the resulting Čech complexes faithfully captures the homotopy type of the underlying data, ensuring topological equivalence between the data and the complex.

\par In the Figures \ref{fig:log-comparison1}, \ref{fig:log-comparison2} and \ref{fig:log-comparison3} one can observe that, for the same synthetic dataset, the distances derived from information theory produce almost the same filtration, generating nearly identical orderings of simplices. When applying our method with the Finsler metric \eqref{lnmetric1}, associated to the distance in \eqref{distanceF_ln}, the resulting filtration also appears very similar.  This outcome is expected, since this particular Finsler metric is isometric to the Euclidean metric and therefore is unlikely to generatea simplicial complex different from that obtained with the Euclidean distance.

\par On the other hand, the Finsler metric introduced in \eqref{alphametric} is not isometric to the Euclidean metric; it depends explicitly on both position and direction. The filtration shown in the Figure \ref{fig:log-comparison4} illustrates a clear departure from the other cases in the figures mentioned in the above paragraph: the metric balls of points near the bottom of the figure remain disjoint until the final stages of the filtration, in contrast to the behavior observed previously. Even in this synthetic dataset, this observation convinces us that the lack of Euclidean isometry and the direction dependence introduce novel features. This suggests that, for real datasets, such metrics may capture previously unrevealed structural information, a question we leave for future work.

\section*{Appendix}
\subsection*{Synthetic Data}
The following points were used in the examples shown in the Figures \eqref{fig:log-comparison1}, \eqref{fig:log-comparison2}, \eqref{fig:log-comparison3} and \eqref{fig:log-comparison4}. Each point lies in the circle in the plane
$\mathbf{n} \cdot \mathbf{x} = 1$, with normal vector $(1,1,1)$, of radius $r=6$, and centered at $(14,15,10)$. These points were synthetically generated to illustrate the differences in the filtration of randomly distributed data. 
They are not derived from any real dataset but were chosen to visualize and compare the behavior of the topological analysis on controlled synthetic examples. The points are:

\begin{enumerate}
    \item $(18.24264069, 15.0, 5.75735931)$
    \item $(16.61833712, 17.27666929, 5.1049936)$
    \item $(14.39420407, 19.03178108, 5.57401485)$
    \item $(12.07976362, 19.86326041, 7.05697597)$
    \item $(10.20522618, 19.58062539, 9.21414843)$
    \item $(9.2000257, 18.2486243, 11.55135001)$
    \item $(9.29444148, 16.17240253, 13.53315599)$
    \item $(10.46684401, 13.82759747, 14.70555852)$
    \item $(12.44864999, 11.7513757, 14.7999743)$
    \item $(14.78585157, 10.41937461, 13.79477382)$
    \item $(16.94302403, 10.13673959, 11.92023638)$
    \item $(18.42598515, 10.96821892, 9.60579593)$
    \item $(18.8950064, 12.72333071, 7.38166288)$
\end{enumerate}

These points were generated to lie on the circle in the plane described above and were used to produce the figures in Section~X of this paper.

%-----------------------------------------------------

\end{document}